\documentclass{amsart}

\theoremstyle{plain}
\newtheorem{thm}{Theorem}[section]
\newtheorem{lem}[thm]{Lemma}

\newtheorem{cor}[thm]{Corollary}
\newtheorem{conj}[thm]{Conjecture}
\theoremstyle{definition}

\theoremstyle{remark}

\newcommand{\C}{\mathbb{C}}

\newcommand{\Z}{\mathbb{Z}}
\newcommand{\N}{\mathbb{N}}
\newcommand{\PP}{\mathbb{P}}

\DeclareMathOperator{\GL}{GL}

\DeclareMathOperator{\tr}{tr}

\DeclareMathOperator{\Hom}{Hom}
\DeclareMathOperator{\ext}{ext}
\DeclareMathOperator{\Ker}{Ker}

\DeclareMathOperator{\Ima}{Im}

\DeclareMathOperator{\Aut}{Aut}
\DeclareMathOperator{\End}{End}
\DeclareMathOperator{\Ext}{Ext}

\DeclareMathOperator{\rank}{rank}
\DeclareMathOperator{\Res}{Res}
\DeclareMathOperator{\conn}{conn}

\newcommand{\OO}{\mathcal{O}}
\newcommand{\E}{\mathbf{E}}

\DeclareMathOperator*{\rightoverleft}{\parbox{2em}{\centerline{$\longrightarrow$}\vskip -6pt\centerline{$\longleftarrow$}}}

\newcommand{\smalltextsc}[1]{#1}

\begin{document}
\title[Indecomposable parabolic bundles]{Indecomposable parabolic bundles and the
existence of matrices in prescribed conjugacy class closures with product equal to the identity}

\dedicatory{Dedicated to Claus Michael Ringel on the occasion of his sixtieth birthday}

\author{William Crawley-Boevey}
\address{Department of Pure Mathematics, University of Leeds, Leeds LS2 9JT, UK}
\email{w.crawley-boevey@leeds.ac.uk}


\begin{abstract}
We study the possible dimension vectors of indecomposable
parabolic bundles on the projective line, and use our answer
to solve the problem of characterizing those collections
of conjugacy classes of $n\times n$ matrices for which one can find matrices
in their closures whose product is equal to the
identity matrix. Both answers depend on the root system of a Kac-Moody Lie algebra.
Our proofs use Ringel's theory of tubular algebras, work of Mihai on the
existence of logarithmic connections, the Riemann-Hilbert correspondence
and an algebraic version, due to Dettweiler and Reiter, of Katz's middle convolution
operation.
\end{abstract}
\maketitle

\section{Introduction}
Given some matrices in $\GL_n(\C)$, if one knows their conjugacy classes,
what can one say about the conjugacy class of their product? This is called
the `recognition problem' by Neto and Silva \cite[\S 3]{NetoS}.
Inverting one of the matrices, one obtains the following more
symmetrical formulation: given conjugacy classes $C_1,\dots,C_k$,
determine whether or not one can solve the equation
\begin{equation}\label{eq1}
A_1 A_2 \dots A_k = 1
\end{equation}
with $A_i\in C_i$. We solve a weaker version of this, determining whether
or not one can solve equation (\ref{eq1}) with matrices $A_i$ in the closures
$\overline{C_i}$ of the conjugacy classes. Of course if the $C_i$ are diagonalizable,
they are already closed, and the recognition problem is solved.

Our solution depends on properties of indecomposable parabolic bundles.
Let $X$ be a connected Riemann surface, $D = (a_1,\dots,a_k)$ a collection
of distinct points of $X$, and $w=(w_1,\dots,w_k)$ a collection of positive integers.
By a \emph{parabolic bundle} on $X$ of \emph{weight type}
$(D,w)$ we mean a collection $\E = (E,E_{ij})$ consisting
of a holomorphic vector bundle $E$ on $X$
and flags of vector subspaces
\[
E_{a_i} = E_{i0} \supseteq E_{i1} \supseteq \dots \supseteq E_{i,w_i-1} \supseteq E_{i,w_i} = 0
\]
of its fibres at the points in $D$.
The parabolic bundles naturally form a category $\mathrm{par}_{D,w}(X)$ in which the
morphisms $\E\to \mathbf{F}$ are the vector bundle homomorphisms $f:E\to F$
with
\[
f_{a_i}(E_{ij}) \subseteq F_{ij}
\]
for all $i,j$. Clearly any points $a_i$ with weighting $w_i=1$ can be
omitted from $D$ without changing the category. Therefore, one may if one wishes
assume that all $w_i\ge 2$. This also shows that in order to specify a weight type,
it suffices to fix the divisor $\sum_{i=1}^k (w_i-1)a_i$.

The category $\mathrm{par}_{D,w}(X)$ has been much studied before.
It is equivalent to the natural category of vector bundles for an
orbifold Riemann surface \cite[Theorem 5.7]{FurSteer}, or,
in case $X$ is the projective line, the category of vector bundles
on a weighted projective line in the sense of Geigle and Lenzing \cite{GL}
(see \cite[Theorem 4.4]{Lenzing}). Analogous categories arise when
one studies one-dimensional smooth stacks or sheaves of torsion-free
$\mathcal{H}$-modules, where $\mathcal{H}$ is a sheaf of classical
hereditary orders on an algebraic curve, see \cite{CI}.

The usual notion of a parabolic bundle due to Mehta and Seshadri \cite{MS,Seshadri}
also associates `weights' in $[0,1)$ to the terms in the flags.
These weights are used to define the notion of semistability for parabolic bundles,
and Mehta and Seshadri relate the moduli space of semistable parabolic bundles on
a compact Riemann surface of genus $g\ge 2$ with
equivalence classes of unitary representations of a Fuchsian group.
Agnihotri and Woodward \cite{AW} and Belkale \cite{Belkale}
independently observed that one can in the same way relate semistable parabolic
bundles on the projective line with unitary representations of the
fundamental group of the punctured projective line. In this way they were
able to solve the recognition problem for conjugacy
classes in $SU(n)$ in terms of quantum Schubert calculus.
In this paper semistability will play no role, so it is convenient to ignore
the weights, and actually, we will need to associate arbitrary complex numbers to the
terms in the flags. Note that the term `quasi-parabolic structure' \cite{MS,Seshadri}
is sometimes used for the flags without the weights.

The basic invariant of a parabolic bundle $\E$ is its \emph{dimension vector}
$\alpha$, which consists of the numbers $\alpha_0 = \rank E$
and $\alpha_{ij}=\dim E_{ij}$ ($1\le i\le k$, $1\le j\le w_i-1$).
Assuming that $X$ is compact,
one also has the degree, $d=\deg E$.
What can one say about the invariants of an \emph{indecomposable}
parabolic bundle?
Clearly the dimension vector must be \emph{strict}, by which we mean that
\begin{equation}
\label{e:strict}
\alpha_0 \ge \alpha_{i1} \ge \alpha_{i2} \ge \dots \ge \alpha_{i,w_i-1} \ge 0
\end{equation}
for all $i$.
If $X$ has genus $g \ge 1$, then there are indecomposable parabolic
bundles of all possible strict dimension vectors and degrees, for
one knows that there is an indecomposable vector bundle
of any rank and degree, and it can be turned into an indecomposable
parabolic bundle using arbitrary flags.
Thus we restrict to the case when $X = \PP^1$, the complex projective
line.

Let $\Gamma_w$ be the star-shaped graph
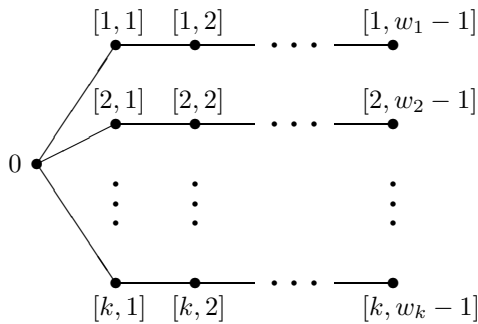
\setlength{\unitlength}{1.5pt}
\[
\begin{picture}(110,80)
\put(10,40){\circle*{2.5}}
\put(30,10){\circle*{2.5}}
\put(30,50){\circle*{2.5}}
\put(30,70){\circle*{2.5}}
\put(50,10){\circle*{2.5}}
\put(50,50){\circle*{2.5}}
\put(50,70){\circle*{2.5}}
\put(100,10){\circle*{2.5}}
\put(100,50){\circle*{2.5}}
\put(100,70){\circle*{2.5}}
\put(10,40){\line(2,-3){20}}
\put(10,40){\line(2,1){20}}
\put(10,40){\line(2,3){20}}
\put(30,10){\line(1,0){35}}
\put(30,50){\line(1,0){35}}
\put(30,70){\line(1,0){35}}
\put(85,10){\line(1,0){15}}
\put(85,50){\line(1,0){15}}
\put(85,70){\line(1,0){15}}
\put(70,10){\circle*{1}}
\put(75,10){\circle*{1}}
\put(80,10){\circle*{1}}
\put(70,50){\circle*{1}}
\put(75,50){\circle*{1}}
\put(80,50){\circle*{1}}
\put(70,70){\circle*{1}}
\put(75,70){\circle*{1}}
\put(80,70){\circle*{1}}
\put(30,25){\circle*{1}}
\put(30,30){\circle*{1}}
\put(30,35){\circle*{1}}
\put(50,25){\circle*{1}}
\put(50,30){\circle*{1}}
\put(50,35){\circle*{1}}
\put(100,25){\circle*{1}}
\put(100,30){\circle*{1}}
\put(100,35){\circle*{1}}
\put(3,38){0}
\put(24,2){$[k,1]$}
\put(24,54){$[2,1]$}
\put(24,74){$[1,1]$}
\put(44,2){$[k,2]$}
\put(44,54){$[2,2]$}
\put(44,74){$[1,2]$}
\put(92,2){$[k,w_k-1]$}
\put(92,54){$[2,w_2-1]$}
\put(92,74){$[1,w_1-1]$}
\end{picture}
\]
with vertex set $I = \{ 0 \} \cup \{ [i,j] : 1\le i\le k, 1\le j\le w_i-1 \}$.
For ease of notation, if $\alpha\in\Z^I$, we write its components as $\alpha_0$ and $\alpha_{ij}$.
In this way one sees that dimension vectors, as defined above, are naturally elements of $\Z^I$.
Associated to $\Gamma_w$ there is a Kac-Moody Lie algebra, and hence a root system, which we
consider as a subset of $\Z^I$. We briefly recall its combinatorial definition. There is a
symmetric bilinear form $(-,-)$ on $\Z^I$ given by
\[
(\epsilon_v,\epsilon_{v'}) = \begin{cases}
2 & \text{(if $v=v'$)} \\
-1 & \text{(if an edge joins $v$ and $v'$)} \\
0 & \text{(otherwise)}
\end{cases}
\]
where $\epsilon_v$ denotes the coordinate vector at a vertex $v\in I$. The
associated quadratic form is $q(\alpha) = \frac12 (\alpha,\alpha)$,
and we also define $p(\alpha)=1-q(\alpha)$. For $v\in I$, the
reflection $s_v : \Z^I \to \Z^I$, is defined by
$s_v(\alpha) = \alpha - (\alpha,\epsilon_v) \epsilon_v$, the \emph{Weyl group}
$W$ is the subgroup of $\Aut(\Z^I)$ generated by the $s_v$, and the
\emph{real roots} are the images under elements of $W$ of the coordinate
vectors. The \emph{fundamental region} consists of the nonzero elements
$\alpha\in\N^I$ which have connected support and which are not made smaller
by any reflection; its closure under the action of $W$ and change of sign is,
by definition, the set of \emph{imaginary roots}. Note that $p(\alpha)=0$ for
real roots, and $p(\alpha) > 0$ for imaginary roots. Recall that any root
is \emph{positive} ($\alpha\in\N^I$), or \emph{negative} ($\alpha\in-\N^I$).

Note that when determining whether or not $\alpha$ is a root,
all that matters is $\alpha_0$ and for each arm $i$, the unordered collection
of numbers
\begin{equation}
\label{e:rootdiffs}
\alpha_0 - \alpha_{i1}, \, \alpha_{i1}-\alpha_{i2}, \, \dots, \, \alpha_{i,w_i-2}-\alpha_{i,w_i-1}, \, \alpha_{i,w_i-1},
\end{equation}
since the reflection at vertex $[i,j]$ has the effect of exchanging the $j$th and $(j+1)$st terms
in this collection.

We say that an element of $\Z^I$ is \emph{strict} if it satisfies the inequalities~(\ref{e:strict}).
Any root with $\alpha_0>0$ must be strict, for if some term in the collection
(\ref{e:rootdiffs}) is negative, it can be moved to the last place, and then $\alpha$ has
both positive and negative components. It follows that
the only positive roots which are not strict are those of the form
\begin{equation}
\label{e:nonstrictroots}
\alpha_0=0,
\quad
\alpha_{ij} = \begin{cases}
1 & \text{(if $i=\ell$ and $r\le j\le s$)} \\
0 & \text{(otherwise),}
\end{cases}
\end{equation}
for some $1\le \ell\le k$ and $1\le r\le s\le w_\ell-1$.

Our first main result is as follows.

\begin{thm}\label{th:mainparroots}
The dimension vector of an indecomposable parabolic bundle on $\PP^1$
of weight type $(D,w)$ is a strict root for $\Gamma_w$.
Conversely, if $\alpha$ is a strict root,
$d\in\Z$, and $d$ and $\alpha$ are coprime, then there is an indecomposable
parabolic bundle $\E$ on $\PP^1$ of weight type $(D,w)$,
with dimension vector $\alpha$ and $\deg E = d$.
Moreover, if $\alpha$ is a real root, this indecomposable parabolic bundle
is unique up to isomorphism.
\end{thm}

Here we say that $d$ and $\alpha$ are \emph{coprime} provided that there is
no integer $\ge 2$ which divides $d$ and all components of $\alpha$.
We conjecture that this coprimality assumption is unnecessary.

Note that in case $\Gamma_w$ is an extended Dynkin diagram,
the corresponding weighted projective line has genus one,
and Lenzing and Meltzer \cite[Theorem 4.6]{LM} have given a complete classification
of indecomposable coherent sheaves on the weighted projective line
in terms of the roots of a quadratic form.

\begin{cor}
The dimension vectors of indecomposable parabolic bundles on $\PP^1$
of weight type $(D,w)$ are exactly the strict roots for $\Gamma_w$.
\end{cor}

We now return to the problem about matrices.
Given a collection of positive integers $w=(w_1,\dots,w_k)$, and a set
$\Omega$, we denote by $\Omega^w$ the set of collections $\xi = (\xi_{ij})$
with $\xi_{ij}\in \Omega$ for $1\le i \le k$ and $1\le j \le w_i$.
To fix conjugacy classes $C_1,\dots,C_k$ in $\GL_n(\C)$ we fix $w$ and
$\xi\in(\C^*)^w$ with
\[
(A_i - \xi_{i1} 1) (A_i - \xi_{i2} 1) \dots (A_i - \xi_{i,w_i} 1) = 0
\]
for $A_i\in C_i$. Clearly, if one wishes one can take $w_i$ to be the
degree of the minimal polynomial of $A_i$, and $\xi_{i1},\dots,\xi_{i,w_i}$
to be its roots, counted with multiplicity. The conjugacy class $C_i$ is
then determined by the ranks of the partial products
\[
\alpha_{ij} = \rank (A_i - \xi_{i1}1)\dots(A_i - \xi_{ij}1)
\]
for $A_i\in C_i$ and $1\le j\le w_i-1$. Setting $\alpha_0 = n$, we obtain
a dimension vector $\alpha$ for the graph $\Gamma_w$.
An obvious necessary condition for a solution to equation (\ref{eq1})
is that
\[
\det A_1 \times \dots \times \det A_k = 1.
\]
For $A_i\in C_i$ this is equivalent to the condition $\xi^{[\alpha]} = 1$,
where we define
\[
\xi^{[\alpha]} = \prod_{i=1}^k \prod_{j=1}^{w_i} \xi_{ij}^{\alpha_{i,j-1} - \alpha_{ij}}
\]
using the convention that $\alpha_{i0} = \alpha_0$ and $\alpha_{i,w_i}=0$ for all $i$.
Our other main result is the following.

\begin{thm}\label{t:clothm}
There is a solution to $A_1\dots A_k = 1$ with $A_i\in\overline{C_i}$
if and only if $\alpha$ can be written as a sum of positive roots for $\Gamma_w$,
say $\alpha = \beta+\gamma+\dots$, with $\xi^{[\beta]} = \xi^{[\gamma]} = \dots = 1$.
\end{thm}

The most basic example, familiar from the hypergeometric equation, is when $k=3$ and
the matrices are $2\times 2$. If $C_i$ has distinct eigenvalues $\lambda_i,\mu_i$,
one can take $w=(2,2,2)$, $\xi_{i1}=\lambda_i$ and $\xi_{i2}=\mu_i$, and
the conjugacy classes are then given by the dimension vector
\[
\begin{picture}(18,40)
\put(0,15){2}
\put(15,0){1}
\put(15,15){1}
\put(15,30){1}
\put(13,5){\line(-1,1){8}}
\put(13,17){\line(-1,0){8}}
\put(13,29){\line(-1,-1){8}}
\end{picture}
\]
for the Dynkin diagram $D_4$. Since this is a positive root, the only requirement
for solvability is $\xi^{[\alpha]}=1$, that is, $\prod_i \lambda_i\mu_i=1$.
On the other hand, if the matrices are $3\times 3$ and $C_i$ is diagonalizable, with $\lambda_i$
having multiplicity 2 and $\mu_i$ having multiplicity 1, the relevant dimension vector is
\[
\begin{picture}(18,40)
\put(0,15){3}
\put(15,0){1}
\put(15,15){1}
\put(15,30){1}
\put(13,5){\line(-1,1){8}}
\put(13,17){\line(-1,0){8}}
\put(13,29){\line(-1,-1){8}}
\end{picture}
\]
which is not a root. Its possible decompositions as a sum of positive roots are
\[
\begin{picture}(108,40)
\put(0,15){1}
\put(15,0){0}
\put(15,15){0}
\put(15,30){1}
\put(13,5){\line(-1,1){8}}
\put(13,17){\line(-1,0){8}}
\put(13,29){\line(-1,-1){8}}
\put(30,15){$+$}
\put(45,15){1}
\put(60,0){0}
\put(60,15){1}
\put(60,30){0}
\put(58,5){\line(-1,1){8}}
\put(58,17){\line(-1,0){8}}
\put(58,29){\line(-1,-1){8}}
\put(75,15){$+$}
\put(90,15){1}
\put(105,0){1}
\put(105,15){0}
\put(105,30){0}
\put(103,5){\line(-1,1){8}}
\put(103,17){\line(-1,0){8}}
\put(103,29){\line(-1,-1){8}}
\end{picture}
\]
\[
\begin{picture}(63,40)
\put(0,15){1}
\put(15,0){0}
\put(15,15){0}
\put(15,30){0}
\put(13,5){\line(-1,1){8}}
\put(13,17){\line(-1,0){8}}
\put(13,29){\line(-1,-1){8}}
\put(30,15){$+$}
\put(45,15){2}
\put(60,0){1}
\put(60,15){1}
\put(60,30){1}
\put(58,5){\line(-1,1){8}}
\put(58,17){\line(-1,0){8}}
\put(58,29){\line(-1,-1){8}}
\end{picture}
\]
and refinements of these, so the condition for solvability is that
$\mu_1\lambda_2\lambda_3=\lambda_1\mu_2\lambda_3=\lambda_1\lambda_2\mu_3=1$
or $\prod_i \lambda_i = \prod_i \lambda_i\mu_i=1$.

Although we have not been able to solve the full recognition problem, we
believe that the approach using indecomposable parabolic bundles will be
of assistance in solving this and related problems. We now make some
remarks about these. A variation of the recognition problem
asks whether or not there is a solution to (\ref{eq1}) which is \emph{irreducible},
meaning that the $A_i$ have no common invariant subspace.
This is called the `Deligne-Simpson problem' by Kostov
\cite{Kostovmon,KostovCR,Kostovaspects,Kostovsurvey}.
Here we have a conjecture.

\begin{conj}
There is an irreducible solution to $A_1\dots A_k = 1$ with $A_i\in C_i$
if and only if $\alpha$ is a positive root, $\xi^{[\alpha]}=1$, and
$p(\alpha) > p(\beta)+p(\gamma)+\dots$ for any nontrivial decomposition
of $\alpha$ as a sum of positive roots $\alpha = \beta+\gamma+\dots$
with $\xi^{[\beta]} = \xi^{[\gamma]} = \dots = 1$.
\end{conj}

A solution to (\ref{eq1}) with $A_i\in C_i$ is said to be \emph{rigid} if
it is the unique such solution, up to simultaneous conjugacy.
Deligne (see \cite[Lemma 6]{Simp}) observed that if there
is an irreducible solution, it is rigid if and only if
\begin{equation}\label{eq:desum}
\sum_{i=1}^k \dim C_i = 2n^2-2.
\end{equation}
In our formulation this becomes the condition $p(\alpha)=0$.
In \cite{Katz}, Katz introduced an operation, called `middle convolution',
which enabled him to give an algorithm for studying rigid irreducible solutions.
In \cite{DR}, Detweiller and Reiter gave a purely algebraic version
of middle convolution, called simply `convolution', and the corresponding algorithm.
In Sections~\ref{s:type}--\ref{s:rigid} we use the methods of Detweiller and Reiter
to prove the rigid case of our conjecture.

\begin{thm}\label{th:rigidcase}
There is a rigid irreducible solution to $A_1\dots A_k = 1$ with $A_i\in C_i$
if and only if $\alpha$ is a positive real root for $\Gamma_w$,
$\xi^{[\alpha]}=1$, and there is no nontrivial decomposition of $\alpha$ as a sum
of positive roots $\alpha = \beta+\gamma+\dots$
with $\xi^{[\beta]} = \xi^{[\gamma]} = \dots = 1$.
\end{thm}

In Section~\ref{s:squid} we study parabolic bundles in
the case when the dimension vector is in the fundamental region.
In Sections~\ref{s:RH}--\ref{s:Weil} we show how to pass between parabolic bundles
and solutions of equation (\ref{eq1}), and then in Section~\ref{s:generic} we
combine all this with the methods used in the rigid case to solve the Deligne-Simpson
problem for generic eigenvalues. The main results, Theorems~\ref{th:mainparroots}
and \ref{t:clothm}, are then deduced from the generic case in Section~\ref{s:mainproofs}.

Finally we remark that instead of equation (\ref{eq1}), one can consider the additive equation
\begin{equation}\label{eq2}
A_1 + A_2 + \dots + A_k = 0.
\end{equation}
We have already solved the corresponding additive Deligne-Simpson problem in \cite{CBad}.
In this case one fixes conjugacy classes
$D_1,\dots,D_k$ in $M_n(\C)$ using an element $\zeta\in\C^w$, and
dimension vector $\alpha$. Defining
\[
\zeta*{[\alpha]} = \sum_{i=1}^k \sum_{j=1}^{w_i} \zeta_{ij}(\alpha_{i,j-1} - \alpha_{ij}),
\]
we showed that there is an irreducible solution to (\ref{eq2}) with $A_i\in D_i$
if and only if $\alpha$ is a positive root for $\Gamma_w$, $\zeta*[\alpha]=0$,
and $p(\alpha) > p(\beta)+p(\gamma)+\dots$ for any nontrivial decomposition
of $\alpha$ as a sum of positive roots $\alpha = \beta+\gamma+\dots$
with $\zeta*[\beta] = \zeta*[\gamma] = \dots = 0$.
Note that the method of \cite{CBad} together with \cite[Theorem 4.4]{CBmm}
and Theorem~\ref{t:conjclotriple} already gives the additive analogue of
Theorem~\ref{t:clothm}. Namely, there is a solution to $A_1+\dots+A_k = 0$ with $A_i\in\overline{D_i}$
if and only if $\alpha$ can be written as a sum of positive roots $\alpha = \beta+\gamma+\dots$
with $\zeta*[\beta] = \zeta*[\gamma] = \dots = 0$.

The crucial Theorem~\ref{t:parcon}, which should be of independent interest,
was proved during a visit to the Center for Advanced Study at the Norwegian
Academy of Science and Letters in September 2001, and I would like to thank
my hosts for their hospitality.
I would also like to thank C.~Gei\ss{} and H.~Lenzing for some invaluable
discussions.

\section{Fixing conjugacy classes and their closures}
\label{s:type}
In Sections \ref{s:type}-\ref{s:squid} we work over an algebraically closed field $K$.

We first deal with a single conjugacy class.
Let $V$ be a vector space of dimension $n$, let $d\ge 1$, and fix a collection
$\xi = (\xi_1,\dots,\xi_d)$ with $\xi_j\in K$.
We say that an endomorphism $A\in \End(V)$ has \emph{type} $\xi$ if
$(A-\xi_1 1)\dots (A-\xi_d 1)=0$.
In this case its \emph{dimension vector}
is the sequence of integers $(n_0,n_1,\dots,n_{d-1})$
where $n_j$ is the rank of the partial product $(A-\xi_1 1)\dots (A-\xi_j 1)$,
so that $n_0=n$. By convention we define $n_d=0$.

Clearly the type and dimension vector only depend on the conjugacy class of $A$.
Now $A$ has type $\xi$ if and only if for all $\lambda\in K$, the Jordan normal form of
$A$ only involves Jordan blocks with eigenvalue $\lambda$ up to size $r_\lambda$, where
\[
r_\lambda = |\{ \ell \mid 1\le \ell\le d, \xi_i = \lambda\}|.
\]
Moreover $n_{j-1} - n_j$ is the number of Jordan blocks of eigenvalue $\xi_j$ of size at
least $m_j$ involved in $A$, where
\[
m_j = |\{ \ell \mid 1\le \ell\le j, \xi_\ell = \xi_j \}|.
\]
It follows that $A\in\End(V)$ of type $\xi$ is determined up to conjugacy by
its dimension vector. Clearly given any $j<\ell$ with $\xi_j = \xi_\ell$
we have $m_j < m_\ell$, so that
\begin{equation}
\label{e:exconjcond}
n_{j-1}-n_{j} \ge n_{\ell-1}-n_{\ell}.
\end{equation}
Conversely any
sequence $(n_0,\dots,n_{d-1})$ of integers with $n=n_0\ge n_1\dots \ge n_{d-1} \ge n_d = 0$
and satisfying (\ref{e:exconjcond}) arises as the
dimension vector of some $A\in\End(V)$ of type $\xi$.

\begin{thm}
\label{t:conjclotriple}
Let $A\in\End(V)$ have type $\xi$ and dimension vector $(n_0,n_1,\dots,n_{d-1})$.
If $B\in\End(V)$, then the following conditions are equivalent.

\textup{\hphantom{ii}(i)}~$B$ is in the closure of the conjugacy class of $A$.

\textup{\hphantom{i}(ii)}~There is a flag of subspaces
\[
V = V_0 \supseteq V_1 \supseteq \dots \supseteq V_d = 0
\]
with $\dim V_j = n_j$ and such that $(B-\xi_j 1)(V_{j-1}) \subseteq V_j$ for all $1\le j\le d$.

\textup{(iii)}~There are vector spaces $V_j$ and linear maps $\phi_j$, $\psi_j$,
\[
V=V_0 \rightoverleft^{\phi_{1}}_{\psi_{1}} V_{1} \rightoverleft^{\phi_{2}}_{\psi_{2}} V_{2}
\rightoverleft^{\phi_{3}}_{\psi_{3}} \dots
\rightoverleft^{\phi_{d}}_{\psi_{d}} V_{d}=0
\]
where $V_j$ has dimension $n_j$, and satisfying
\[
\begin{split}
B - \psi_{1} \phi_{1}
&= \xi_{1} 1,
\\
\phi_{j} \psi_{j} - \psi_{j+1}\phi_{j+1}
&= (\xi_{j+1}-\xi_{j}) 1,
\qquad (1\le j < d).
\end{split}
\]
\end{thm}

In case $K$ has characteristic 0, the equivalence
of (i) and (iii) follows from \cite[Lemma 9.1]{CBnorm}.

\begin{proof}
(i)$\Rightarrow$(ii).
Let $F$ be the corresponding flag variety, and let $Z$ be the closed
subset of $F\times\End(V)$ consisting of all flags $V_j$ and endomorphisms
$B$ with $(B-\xi_j 1)(V_{j-1}) \subseteq V_j$ for all $1\le j\le d$.
By definition the set of $B$ which satisfy condition (ii) is the image $I$
of $Z$ under the projection to $\End(V)$.
Now $I$ contains $A$, since one can take $V_j = \Ima(A-\xi_1 1)\dots (A-\xi_j 1)$,
it is stable under the conjugation action of $\GL(V)$, and
it is closed since $F$ is projective. Thus it contains the closure
of the conjugacy class of $A$.

(ii)$\Rightarrow$(iii).
Take $\phi_j$ to be the restriction of $B-\xi_j 1$ to $V_{j-1}$ and $\psi_j$ to be the inclusion.

(iii)$\Rightarrow$(i).
We prove this by induction on $d$. If $d=1$ it is trivial, so suppose that $d>1$.
Identify $V_1$ with $\Ima(A-\xi_1 1)$, and let $A_1$ be the restriction of $A$ to $V_1$.
Observe that $A_1$ has type $(\xi_2,\dots,\xi_d)$ and dimension vector $(n_1,\dots,n_d)$.
Now letting $B_1 = \phi_1 \psi_1 + \xi_1 1$, we have
\[
\begin{split}
B_1 - \psi_2 \phi_2
&= \xi_2 1,
\\
\phi_{j} \psi_{j} - \psi_{j+1}\phi_{j+1}
&= (\xi_{j+1}-\xi_{j}) 1,
\qquad (2\le j < d)
\end{split}
\]
so by the inductive hypothesis $B_1$ is in the closure of the conjugacy class in $\End(V_1)$
containing $A_1$.

Recall that the Gerstenhaber-Hesselink Theorem \cite[Theorem 1.7]{Gerstenhaber}
says that if $A,B\in\End(V)$, then $B$ is in the closure of the conjugacy class of $A$ if and only if
\[
\rank (B-\lambda 1)^m \le \rank (A-\lambda 1)^m
\]
for all $\lambda\in K$ and $m\ge 1$. Applying this to $A_1$ and $B_1$, we have
\[
\rank  (B_1-\lambda 1)^m \le \rank (A_1-\lambda 1)^m = \rank (A-\xi_1 1)(A-\lambda 1)^m.
\]
Now if $\lambda=\xi_1$ we have
\[
(B-\lambda1)^{m+1} = (\psi_1 \phi_1)^{m+1} = \psi_1 (B_1 - \lambda 1)^m \phi_1
\]
so that
\[
\rank (B-\lambda1)^{m+1} \le \rank (B_1 - \lambda 1)^m \le \rank (A-\lambda 1)^{m+1}.
\]
On the other hand, if $\lambda\neq \xi_1$, then
\[
\rank(\psi_1\phi_1 + (\xi_1-\lambda)1)^m = n_0 - n_1 + \rank(\phi_1\psi_1 +(\xi_1-\lambda)1)^m
\]
by Lemma~\ref{l:mapsback} below, so
$\rank (B-\lambda 1)^m = n_0 - n_1 + \rank (B_1 - \lambda 1)^m$.
However, we also have
\[
\rank (A_1 - \lambda)^m = \rank (A-\xi_1 1)(A-\lambda 1)^m = n_1 - n_0 + \rank(A-\lambda 1)^m
\]
using that
\[
\Ker (A-\xi_1 1)(A-\lambda1)^m  = \Ker(A-\xi_1 1)\oplus \Ker(A-\lambda1)^m,
\]
and it follows that $\rank (B-\lambda1)^m \le \rank (A-\lambda1)^m$.
Since this inequality holds for all $\lambda$ and $m$, $B$ is in the closure of the conjugacy class of $A$
by the Gerstenhaber-Hesselink Theorem.
\end{proof}

\begin{lem}
\label{l:mapsback}
If $\phi:V\to W$ and $\psi:W\to V$ are linear maps, then
\[
\dim V - \rank (\psi\phi + \mu 1)^m = \dim W - \rank (\phi\psi + \mu 1)^m
\]
for any $0\neq\mu\in K$ and $m\ge 0$.
\end{lem}

\begin{proof}
From the formula $\psi(\phi\psi + \mu 1) = (\psi\phi+\mu 1) \psi$
one deduces that
\[
\psi(\phi\psi + \mu 1)^m = (\psi\phi+\mu 1)^m \psi,
\]
so that $\psi$ induces a map from $\Ker(\phi\psi + \mu 1)^m$ into $\Ker (\psi\phi+\mu 1)^m$.
This map is injective, since if $\psi(x)=0$, then $(\phi\psi + \mu 1)^m (x)=\mu^m x$, so $x=0$.
Thus
\[
\dim \Ker(\phi\psi + \mu 1)^m \le \dim\Ker (\psi\phi+\mu 1)^m.
\]
The reverse inequality holds by symmetry, and the result follows.
\end{proof}

We consider the following operation on a sequence $\mathbf{n} = (n_0,n_1,\dots,n_{d-1})$.
Suppose that $1\le r\le s\le d-1$, that $\xi_r = \xi_{s+1}$, and that
the sequence
\[
\mathbf{n'} = (n_0,\dots,n_{r-1},n_r-1,n_{r+1}-1,\dots n_s-1,n_{s+1},\dots)
\]
still has all terms non-negative. In this case we say that
the sequence $\mathbf{n'}$ is obtained from $\mathbf{n}$ by \emph{reduction} with respect to $r,s$ (and $\xi$).

\begin{lem}
\label{l:oddconjclo}
Suppose that $A\in\End(V)$ has type $\xi$ and dimension vector $\mathbf{n} = (n_0,n_1,\dots,n_{d-1})$.
Suppose that $\mathbf{m} = (m_0,m_1,\dots,m_{d-1})$ is obtained from $\mathbf{n}$
by a finite number of reductions.
Let $B\in\End(V)$.
If there are vector spaces $V^0_j$ and linear maps $\phi^0_j$, $\psi^0_j$,
\[
V=V^0_0 \rightoverleft^{\phi^0_{1}}_{\psi^0_{1}} V^0_{1} \rightoverleft^{\phi^0_{2}}_{\psi^0_{2}} V^0_{2}
\rightoverleft^{\phi^0_{3}}_{\psi^0_{3}} \dots
\rightoverleft^{\phi^0_{d}}_{\psi^0_{d}} V^0_{d}=0
\]
where $V^0_j$ has dimension $m_j$, and satisfying
\[
\begin{split}
B - \psi^0_{1} \phi^0_{1}
&= \xi_{1} 1,
\\
\phi^0_{j} \psi^0_{j} - \psi^0_{j+1}\phi^0_{j+1}
&= (\xi_{j+1}-\xi_{j}) 1,
\qquad (1\le j < d),
\end{split}
\]
then $B$ is in the closure of the conjugacy class of $A$.
\end{lem}

\begin{proof}
Suppose that $\mathbf{m}$ is obtained from $\mathbf{n}$ by reducing with respect to $r_i$ and $s_i$
for $1\le i\le p$. For each $1\le i\le p$, define vector spaces
\[
V^i_j = \begin{cases}
K & \text{($r_i \le j \le s_i$)} \\
0 & \text{(otherwise)}
\end{cases}
\]
and linear maps
\[
0=V^i_0 \rightoverleft^{\phi^i_{1}}_{\psi^i_{1}} V^i_{1} \rightoverleft^{\phi^i_{2}}_{\psi^i_{2}} V^i_{2}
\rightoverleft^{\phi^i_{3}}_{\psi^i_{3}} \dots
\rightoverleft^{\phi^i_{d}}_{\psi^i_{d}} V^0_{d}=0
\]
by $\phi^i_j = (\xi_j - \xi_{r_i}) 1$ and $\psi^i_j = 1$ for $r_i< j\le s_i$ and
$\phi^i_j = 0$ and $\psi^i_j=0$ otherwise.
Bearing in mind that many of the spaces are zero, and
using the fact that $\xi_{r_i} = \xi_{s_i+1}$, one easily sees that
\[
\begin{split}
 - \psi^i_{1} \phi^i_{1}
&= \xi_{1} 1,
\\
\phi^i_{j} \psi^i_{j} - \psi^i_{j+1}\phi^i_{j+1}
&= (\xi_{j+1}-\xi_{j}) 1,
\qquad (1\le j < d).
\end{split}
\]
Now the vector spaces
\[
V_j = \bigoplus_{i=0}^p V^i_p
\]
have dimension $n_j$, and the maps $\phi_j$ and $\psi_j$ given by the direct sums of the $\phi^i_j$ and
the $\psi^i_j$ satisfy condition (iii) in Theorem~\ref{t:conjclotriple}. Thus $B$ is in
the closure of the conjugacy class of $A$.
\end{proof}

Now we introduce notation for dealing with collections of conjugacy classes.
Let $k\ge 0$ and let $V_1,\dots,V_k$ be $n$-dimensional vector spaces.
We consider $k$-tuples $A_1,\dots,A_k$ with $A_i\in\End(V_i)$.
Allowing the possibility that each $V_i=K^n$, this may just be a
$k$-tuple of $n\times n$-matrices.

If $w=(w_1,\dots,w_k)$ is a collection of positive integers,
and $\xi\in K^w$, we say that a $k$-tuple $A_1,\dots,A_k$ has \emph{type} $\xi$
provided that $\prod_{j=1}^{w_i}(A_i-\xi_{ij}1)=0$ for all $i$. If so, then as in the
introduction we define the \emph{dimension vector} of the $k$-tuple
(with respect to $\xi$) to be the element $\alpha\in\Z^I$, where $I$ is the
vertex set of the graph $\Gamma_w$, given by $\alpha_0 = n$ and
$\alpha_{ij}=\rank \prod_{\ell=1}^j (A_i - \xi_{i\ell}1)$.

Note that not all dimension vectors $\alpha\in\N^I$ arise this way.
By (\ref{e:exconjcond}), the $\alpha$ which can arise are the strict
elements with
\begin{equation}
\label{e:conjdimveccond}
\alpha_{i,j-1}-\alpha_{ij} \ge \alpha_{i,\ell-1}-\alpha_{i\ell}
\end{equation}
for all $i$ and $j<\ell$ with $\xi_{ij}=\xi_{i\ell}$.

Clearly the type and dimension vector of a $k$-tuple are unchanged by conjugation
of any of the terms. Thus one can speak about the type of a $k$-tuple
of conjugacy classes in $\End(V_i)$ or $M_n(K)$.

\section{Convolution}
In order to study local systems on a punctured Riemann sphere,
Katz \cite{Katz} defined a `middle convolution' operator.
A purely algebraic version, called simply `convolution', was found by
Dettweiler and Reiter \cite{DR}. In this section we describe the
effect of this operator, using our notation for conjugacy classes.
(An alternative algebraic version of middle convolution was found by V\"olklein \cite{V}.)

Let $G$ be the group $\langle g_1,\dots,g_k : g_1\dots g_k=1 \rangle$.
If $\xi\in(K^*)^w$, we say that a representation $\rho:G\to\GL_n(K)$ is of
type $\xi$ if the $k$-tuple $\rho(g_1),\dots,\rho(g_k)$ is of type $\xi$, and we
define the dimension vector correspondingly. Note that the dimension vector is
additive on direct sums, but not necessarily on short exact sequences.
The following fact was mentioned in the introduction.

\begin{lem}\label{l:repdetcond}
If $\rho$ is a representation of $G$ of type $\xi$ and
dimension vector $\alpha$, then $\xi^{[\alpha]}=1$.
\end{lem}

We say that an irreducible representation of $G$ is \emph{$\xi$-collapsing}
if it is of type $\xi$, 1-dimensional, and at least $k-1$ of the $g_i$ act as $\xi_{i1} 1$.
Thus its dimension vector $\alpha$ is either $\epsilon_0$, or there are $1\le \ell\le k$
and $1\le r\le w_\ell-1$ such that
\begin{equation}
\label{e:colldimv}
\alpha_0=1,
\quad
\alpha_{ij} = \begin{cases}
1 & \text{($i=\ell$ and $j\le r$)} \\
0 & \text{(otherwise).}
\end{cases}
\end{equation}
Observe that $s_0(\alpha)$ is either negative, or not strict, hence the name `collapsing'.
We say that a representation is \emph{$\xi$-noncollapsing} if it is of type $\xi$,
and it has no subrepresentation or quotient which is $\xi$-collapsing.

Define $r_0':(K^*)^w\to (K^*)^w$ by
\[
r_0'(\xi)_{ij} = \begin{cases}
1/\xi_{i1} & \text{(if $j=1$)} \\
\frac{\xi_{ij}}{\xi^2_{i1}} \prod_{s=1}^k \xi_{s1} & \text{(if $j>1$).}
\end{cases}
\]
Observe that $r'_0(\xi)^{[s_0(\alpha)]}=\xi^{[\alpha]}$ for any $\alpha$.

Our formulation of convolution is as follows.

\begin{thm}[Dettweiler and Reiter]\label{t:conv}
Given $\xi\in(K^*)^w$, if $\xi^{[\epsilon_0]}\neq 1$, then
there is an equivalence $R_0$ from the category of $\xi$-noncollapsing
representations of $G$ to the category of $r_0'(\xi)$-noncollapsing representations
of $G$. It acts as on dimension vectors as the reflection $s_0$. Moreover $R_0(\rho)$
is irreducible if and only if $\rho$ is irreducible.
\end{thm}

\begin{proof}
Let $\lambda = \xi^{[\epsilon_0]} = \prod_i \xi_{i1} \neq 1$.
For simplicity we write $\xi'$ for $r_0'(\xi)$.
Given a representation $\rho:G\to\GL(V)$, define $A_1,\dots,A_k$
by
\[
A_{k+1-i} = \frac{1}{\xi_{i1}}\, \rho(g_i) \in \GL(V),
\]
so that $\lambda A_k \dots A_2 A_1 = 1$.
As in \cite{DR}, define $G_i\in \GL(V^k)$ by the block matrix
\[
G_i = \begin{pmatrix}
1 & \dots & 0 & 0 & 0 & \dots & 0 \\
\vdots & \ddots & \vdots & \vdots & \vdots & & \vdots \\
0 & \dots & 1 & 0 & 0 & \dots & 0 \\
A_1-1 & \dots & A_{i-1}-1 & \lambda A_i & \lambda(A_{i+1}-1) & \dots & \lambda(A_k-1) \\
0 & \dots & 0 & 0 & 1 & \dots & 0 \\
\vdots & & \vdots & \vdots & \vdots & \ddots & \vdots \\
0 & \dots & 0 & 0 & 0 & \dots & 1
\end{pmatrix}.
\]
Let also $D\in \GL(V^k)$ be given by
\[
D =\begin{pmatrix}
A_1 & 0 & \dots & 0 \\
0 & A_2 & \dots & 0 \\
\vdots & \vdots & \ddots & \vdots \\
0 & 0 & \dots & A_k
\end{pmatrix}.
\]
Define subspaces of $V^k$,
\[
\mathcal{K}_V = \Ker(D-1),
\quad\text{and}\quad
\mathcal{L}_V = \Ker(G_1-1)\cap \dots \cap \Ker(G_k-1).
\]
These are invariant subspaces for the $G_i$.
We define $R_0(\rho)$ to be the representation $\rho':\langle g_1,\dots,g_k\rangle \to\GL(V')$,
where
\[
V' = V^k / (\mathcal{K}_V+\mathcal{L}_V)
\quad\text{and}\quad
\rho'(g_{k+1-i}) = \frac{1}{\xi_{k+1-i,1}} \, \overline{G_i}
\]
where $\overline{G_i}$ denotes the endomorphism of $V'$ induced by $G_i$.
This defines $R_0$ on objects, but it is clear how to extend this
to a functor, see \cite[Proposition 2.6]{DR}.

In \cite[\S 3]{DR}, Dettweiler and Reiter consider the two conditions
\begin{align*}
&(*)  & \bigcap_{j\neq i} \Ker(A_j - 1) \cap \Ker(\tau A_i-1) &= 0, & (i&=1,\dots,k, \; \tau\in K^*), \text{ and}
\\
&(**) & \sum_{j\neq i} \Ima(A_j-1) + \Ima(\tau A_i-1)         &= V, & (i&=1,\dots,k, \; \tau\in K^*).
\end{align*}
Observe that the $A_i$ satisfy condition (*) if and only if $V$ has no
nonzero subspace on which at least $k-1$ of the
$g_i$ act as multiplication by $\xi_{i1}$, and the remaining $g_i$ acts
as multiplication by some scalar. Thus the $A_i$ satisfy condition (*)
if and only if $V$ has no $\xi$-collapsing subrepresentation.
Dually,
it is easy to see that the $A_i$ satisfy (**) if and only if $V$ has no
$\xi$-collapsing quotient representation. Thus $V$ is $\xi$-noncollapsing if and only if
the $A_i$ satisfy (*) and (**).

Now suppose that $V$ is $\xi$-noncollapsing. Clearly
\[
U_\infty = \Ima(A_1-1) + \dots + \Ima(A_k-1)
\]
is an invariant subspace for the $A_i$. All
the $A_i$ act trivially on the quotient $V/U_\infty$, but
$\lambda A_k\dots A_1=1$, and $\lambda\neq 1$, and hence $U_\infty=V$.
Thus \cite[Lemma 4.1(b)]{DR} implies that $\overline{G_k} \dots \overline{G_1} = \lambda 1$,
so that $\rho'(g_1)\dots \rho'(g_k) = 1$, and $R_0(V)$ is actually a representation of $G$.

We show that $R_0(V)$ has type $\xi'$ and that if $V$ has dimension $\alpha$, then $V'$ has
dimension $\alpha'=s_0(\alpha)$.
First, by \cite[Lemma 2.7]{DR}, $\mathcal{L}_V$ is isomorphic to $V$, as a vector
space. Also
\[
\rho(g_i) - \xi_{i1} 1 = \xi_{i1} (A_{k+1-i}-1)
\]
so that
\[
\dim V' = \sum_{i=1}^k \rank(\rho(g_i) - \xi_{i1} 1) - \dim \mathcal{L} = \sum \alpha_{i1} - \alpha_0 = \alpha'_0.
\]
Also,
\[
\rho'(g_i)-\xi'_{i1}1 = \frac{1}{\xi_{i1}} (\overline{G_{k+1-i}}-1),
\]
and
\[
\rank(\overline{G_{k+1-i}}-1) = \rank(A_{k+1-i}-1)
\]
by \cite[Lemma 4.2(a)]{DR}, so that $\alpha'_{i1} = \alpha_{i1}$.
Now we have
\begin{align*}
&\rank \prod_{j=1}^\ell(\rho(g_i)-\xi_{ij}1)
= \rank \prod_{j=2}^\ell(\rho(g_i)-\xi_{ij}1)|_{\Ima (\rho(g_i)-\xi_{i1}1)}
\\
&\quad= \rank \prod_{j=2}^\ell(A_{k+1-i}-\frac{\xi_{ij}}{\xi_{i1}} \, 1)
|_{\Ima (A_{k+1-i}-1)}
\\
&\quad=\rank \prod_{j=2}^\ell
(\overline{G_{k+1-i}}-\frac{\xi_{ij}\lambda}{\xi_{i1}} \, 1)
|_{\Ima (\overline{G_{k+1-i}}-1)}
\text{ by \cite[Lemma 4.1(a)]{DR}}
\\
&\quad= \rank \prod_{j=2}^\ell
(\rho'(g_i)-\frac{\xi_{ij}\lambda}{\xi_{i1}^2} \, 1)
|_{\Ima (\rho'(g_i)-\frac{1}{\xi_{i1}} \, 1)}
\\
&\quad= \rank \prod_{j=1}^\ell (\rho'(g_i)-\xi'_{ij} 1).
\end{align*}
It follows that $R_0(V)$ is of type $\xi'$ and $\alpha'_{ij} = \alpha_{ij}$
for $j>1$. Thus $\alpha' = s_0(\alpha)$, as claimed.

Note also that $R_0(V)$ is $\xi'$-noncollapsing if and only if the $\overline{G_i}$
satisfy (*) and (**), so \cite[Proposition]{DR} shows that $R_0$ sends $\xi$-noncollapsing
representations of $G$ to $\xi'$-noncollapsing representations of $G$.

If we temporarily denote this functor by $R_0^\xi$, then there is also
a functor $R_0^{\xi'}$ sending $\xi'$-noncollapsing
representations of $G$ to $\xi$-noncollapsing representations of $G$
(because clearly $r_0'(\xi')=\xi$). By \cite[Theorem 3.5 and Proposition 3.2]{DR},
these functors are inverse equivalences.
Now \cite[Corollary 3.6]{DR} completes the proof.
\end{proof}

\section{Rigid case of the Deligne-Simpson Problem}\label{s:rigid}
Let $G = \langle g_1,\dots,g_k : g_1\dots g_k=1 \rangle$.
If $\xi\in (K^*)^w$ and $\alpha\in\N^I$, we say that a representation
of $G$ of type $\xi$ and dimension vector $\alpha$ is \emph{rigid}
if, up to isomorphism, it is the only representation of type $\xi$ and dimension vector $\alpha$.

Given $\xi\in (K^*)^w$, we write $S_\xi$ for the set of
strict real roots $\alpha$ for $\Gamma_w$
with $\xi^{[\alpha]}=1$, and such that
there is no nontrivial decomposition of $\alpha$ as a sum of positive roots
$\alpha = \beta+\gamma+\dots$
with $\xi^{[\beta]} = \xi^{[\gamma]} = \dots = 1$.
We show that $S_\xi$ is the set of dimension vectors of rigid representations
of $G$ of type $\xi$.

\begin{lem}\label{l:decond}
An irreducible representation of $G$ of type $\xi$ and dimension vector $\alpha$
is rigid if and only if $p(\alpha)=0$.
\end{lem}

\begin{proof}
We mentioned this in the introduction.
The equivalence of equation (\ref{eq:desum}) and $p(\alpha)=0$
is easy. Thus the lemma follows from Deligne's observation in \cite{Simp}.
Note that this holds for an arbitrary algebraically closed field---see \cite{SV}.
\end{proof}

\begin{lem}[Scott]
\label{l:scott}
If an irreducible representation of $G$ has type $\xi$, dimension $\alpha\neq \epsilon_0$,
and $\xi^{[\epsilon_0]}=1$, then $2\alpha_0 \le \sum_{i=1}^k \alpha_{i1}$, so that
$s_0(\alpha) \ge \alpha$.
\end{lem}

\begin{proof}
Let the representation be $\rho:G\to\GL_n(K)$, with $n=\alpha_0$.
The condition on $\xi$ implies that one obtains another
irreducible representation $\sigma:G\to\GL_n(K)$ by defining
$\sigma(g_i)=\frac{1}{\xi_{i1}} \rho(g_i)$.
Now the condition that $\alpha\neq\epsilon_0$ ensures that
$\sigma$ is not the 1-dimensional trivial representation.
Then, since it is irreducible, a result of Scott \cite[Theorem 1]{Scott} gives
\[
2n \le \sum_{i=1}^k \rank(\sigma(g_i)-1) = \sum_{i=1}^k \rank(\rho(g_i)-\xi_{i1} 1) = \sum_{i=1}^k \alpha_{i1},
\]
as required.
\end{proof}

The convolution functor acts on dimension vectors as the reflection $s_0$
corresponding to the vertex 0 in $\Gamma_w$. We show that reflections at other vertices
can be obtained by permuting the components of $\xi\in (K^*)^w$.
Suppose $1\le i\le k$ and $1\le j\le w_i-1$, and let $v=[i,j]$ be the
corresponding vertex in $\Gamma_w$. Define
$r_v':(K^*)^w \to (K^*)^w$ by
\[
r_v'(\xi)_{pq} = \begin{cases}
\xi_{ij} & \text{(if $p=i$ and $q=j+1$)} \\
\xi_{i,j+1} & \text{(if $p=i$ and $q=j$)} \\
\xi_{pq} & \text{(otherwise).}
\end{cases}
\]
Note that $r'_v(\xi)^{[s_v(\alpha)]}=\xi^{[\alpha]}$ for any $\alpha$.

\begin{lem}\label{l:armreflrep}
Let $v=[i,j]$ and let $\alpha\in\N^I$.

\textup{\hphantom{i}(i)}~If there is a representation of type $\xi$ and dimension vector $\alpha$, and
$s_v(\alpha) < \alpha$, then $\xi_{ij}\neq \xi_{i,j+1}$, that is, $\xi^{[\epsilon_v]}\neq 1$.

\textup{(ii)}~If $\xi^{[\epsilon_v]}\neq 1$ then representations of type $\xi$ with dimension vector $\alpha$
are exactly the same as representations of type $r_v'(\xi)$ with dimension vector $s_v(\alpha)$.
\end{lem}

\begin{proof}
(i) This follows from equation (\ref{e:conjdimveccond}).

(ii) If $\xi^{[\epsilon_v]}\neq 1$, then $\xi_{ij}\neq \xi_{i,j+1}$.
Since $r'_v(\xi)$ is obtained from $\xi$ by exchanging $\xi_{ij}$ and $\xi_{i,j+1}$,
it follows from the definition that a representation $\rho$ has type $\xi$
if and only if it has type $r'_v(\xi)$.
Suppose that $\rho$ has dimension vector $\alpha$ with respect to $\xi$ and $\alpha'$ with
respect to $r_v'(\xi)$. Clearly we have $\alpha_{i\ell} = \alpha'_{i\ell}$ for $\ell\neq j$,
since the products $(\rho(g_i) - \xi_{i1}1)\dots (\rho(g_i) - \xi_{i\ell}1)$,
and $(\rho(g_i) - r_v'(\xi)_{i1}1)\dots (\rho(g_i) - r_v'(\xi)_{i\ell}1)$
are the same, except possibly that two terms are exchanged.
To relate $\alpha_{ij}$ and $\alpha'_{ij}$, let
$W = \Ima (\rho(g_i) - \xi_{i1}1)\dots (\rho(g_i) - \xi_{i,j-1}1)$,
so $\dim W = \alpha_{i,j-1}$.
Then
\[
\alpha_{ij} = \dim (\rho(g_i) - \xi_{ij}1)(W) = \alpha_{i,j-1} - x,
\]
where $x$ is the multiplicity of $\xi_{ij}$ as an eigenvalue of $\rho(g_i)$ on $W$.
Moreover, since $\xi_{ij}\neq \xi_{i,j+1}$, we have
\[
\alpha_{i,j+1} = \dim (\rho(g_i) - \xi_{ij}1)(\rho(g_i) - \xi_{i,j+1}1)(W) = \alpha_{i,j-1} - x- y,
\]
where $y$ is the multiplicity of $\xi_{i,j+1}$ as an eigenvalue of $\rho(g_i)$ on $W$.
On the other hand
\[
\alpha'_{ij} = \dim (\rho(g_i) - \xi_{i,j+1}1)(W) = \alpha_{i,j-1} - y
= \alpha_{i,j-1} + \alpha_{i,j+1} - \alpha_{ij}.
\]
Thus $\alpha' = s_v(\alpha)$.
\end{proof}

\begin{lem}\label{l:reflS}
Let $\xi\in(K^*)^w$ and $\alpha\in \N^I$.
Let $v$ be a vertex in $\Gamma_w$.

\textup{\hphantom{i}(i)}~If $\alpha\in S_\xi$, $\alpha\neq \epsilon_0$ and $s_v(\alpha) < \alpha$, then $\xi^{[\epsilon_v]} \neq 1$.

\textup{(ii)}~If $\alpha$ and $s_v(\alpha)$ are strict and $\xi^{[\epsilon_v]} \neq 1$,
then $\alpha\in S_\xi\Leftrightarrow s_v(\alpha)\in S_{r_v'(\xi)}$.
\end{lem}

\begin{proof}
(i) Since $\alpha$ is strict and $\alpha\neq \epsilon_0$, we have $\alpha\neq \epsilon_v$,
so that $s_v(\alpha)$ is a positive root.
But now there is a decomposition
\[
\alpha = s_v(\alpha) + \epsilon_v + \dots + \epsilon_v.
\]
Since $\xi^{[\alpha]}=1$, if $\xi^{[\epsilon_v]}=1$, then also $\xi^{[s_v(\alpha)]}=1$,
and this decomposition contradicts the fact that $\alpha\in S_\xi$.

(ii) By symmetry it suffices to prove one implication. Applying $s_v$ to any decomposition
of $\alpha$ as a sum of positive roots, one obtains a decomposition of $s_v(\alpha)$ as
a sum of roots. The only difficulty is to ensure that they are all positive. But the
only positive root which changes sign is $\epsilon_v$, and this is ruled out since $\xi^{[\epsilon_v]}\neq 1$.
\end{proof}

Theorem~\ref{th:rigidcase} can now be formulated as follows.

\begin{thm}\label{t:rigidds}
If $\xi\in (K^*)^w$, then there is a rigid irreducible representation of $G$ of
type $\xi$ and dimension $\alpha$ if and only if $\alpha\in S_\xi$.
\end{thm}

\begin{proof}
We prove this for all $\xi$ by induction on $\alpha$.
Either condition requires $\alpha$ to be strict, so the smallest
possible dimension vector is the coordinate vector $\epsilon_0$.
For the equivalence in this case, observe that
$\xi^{[\alpha]}=\prod_{i=1}^k \xi_{i1}$, and if this
is equal to 1, then the representation $\rho(g_i)=\xi_{i1} 1$
is clearly irreducible and rigid.
Thus suppose that $\alpha\neq \epsilon_0$.

By Lemma~\ref{l:decond}, either the existence of a rigid
irreducible representation, or the assumption that $\alpha\in S_\xi$
implies that $p(\alpha)=0$. Then
\[
\sum_{v\in I} \alpha_v(\alpha,\epsilon_v) = (\alpha,\alpha) = 2,
\]
so there must be a vertex $v$ with $(\alpha,\epsilon_v) > 0$, and therefore $s_v(\alpha)<\alpha$.

If there is a vertex of the form $v=[i,j]$ with $s_v(\alpha) < \alpha$,
then Lemmas~\ref{l:armreflrep} and \ref{l:reflS},
show that either condition ensures that $\xi^{[\epsilon_v]}\neq 1$, and they then reduce
the problem to the corresponding one for $r'_v(\xi)$ and $s_v(\alpha)$, so one can use induction.

Thus suppose that the only vertex $v$ with $s_v(\alpha)<\alpha$ is $v=0$.
In particular $\alpha$ is not of the form (\ref{e:colldimv}),
which is reduced by the reflection at $v=[\ell,r]$.

If there is a rigid irreducible representation of $G$ of
type $\xi$ and dimension vector $\alpha$, then
$\xi^{[\epsilon_0]}\neq 1$ by Lemma~\ref{l:scott}.
Moreover, since $\alpha$ is not $\epsilon_0$ or of the form (\ref{e:colldimv}),
the representation is not $\xi$-collapsing, so since it
is irreducible, it is $\xi$-noncollapsing.
Thus by convolution, Theorem~\ref{t:conv}, there is an irreducible
representation of type $r_0'(\xi)$ and dimension vector $s_0(\alpha)$.
It is rigid by Lemma~\ref{l:decond}, since $p(s_0(\alpha)) = p(\alpha)$.
Then by induction $s_0(\alpha)\in S_{r_0'(\xi)}$, and so $\alpha\in S_\xi$
by Lemma~\ref{l:reflS}.

Conversely, suppose that $\alpha\in S_\xi$. Since $\alpha$ is a positive root
and $\alpha\neq\epsilon_0$, $s_0(\alpha)$ is a positive root. Moreover $s_0(\alpha)$ is strict,
for otherwise it must be of the form (\ref{e:nonstrictroots}), and hence
$\alpha$ must be of the form (\ref{e:colldimv}).
Thus by Lemma~\ref{l:reflS}, $\xi^{[\epsilon_0]}\neq 1$ and
$s_0(\alpha)\in S_{r_0'(\xi)}$. Now by induction there is
a rigid irreducible representation of $G$ of
type $r'_0(\xi)$ and dimension vector $s_0(\alpha)$.
This representation is not $r'_0(\xi)$-collapsing,
since $\alpha = s_0(s_0(\alpha))$ is strict, so since it is
irreducible it is $r'_0(\xi)$-noncollapsing. Thus by
Theorem~\ref{t:conv} there is an irreducible representation
of $G$ of type $\xi$ and dimension vector $\alpha$.
It is rigid by Lemma~\ref{l:decond}.
\end{proof}

\section{Squids and the fundamental region}\label{s:squid}
To parametrize families of parabolic bundles on $\PP^1$
it is possible to use a quot-scheme construction, but
it is equivalent to study a space of
representations of a certain algebra $S_{D,w}$,
called a \emph{squid}. By combining general position arguments with
known results about the representation theory of the squid,
we will show the existence of indecomposable parabolic bundles in some cases.

Given a collection $D = (a_1,\dots,a_k)$ of distinct points of $\PP^1$,
and a collection $w=(w_1,\dots,w_k)$ of positive integers, we define
$S_{D,w}$ to be the finite-dimensional associative algebra
given by the quiver $Q_w$
\setlength{\unitlength}{1.5pt}
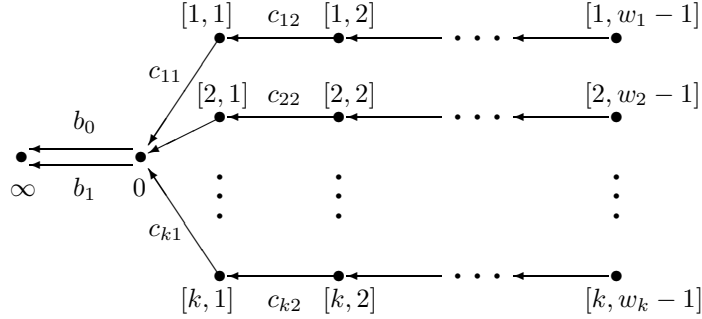
\[
\begin{picture}(110,80)
\put(-20,40){\circle*{2.5}}
\put(10,40){\circle*{2.5}}
\put(30,10){\circle*{2.5}}
\put(30,50){\circle*{2.5}}
\put(30,70){\circle*{2.5}}
\put(60,10){\circle*{2.5}}
\put(60,50){\circle*{2.5}}
\put(60,70){\circle*{2.5}}
\put(130,10){\circle*{2.5}}
\put(130,50){\circle*{2.5}}
\put(130,70){\circle*{2.5}}
\put(8,38){\vector(-1,0){26}}
\put(8,42){\vector(-1,0){26}}
\put(30,70){\vector(-2,-3){18}}
\put(30,50){\vector(-2,-1){18}}
\put(30,10){\vector(-2,3){18}}
\put(60,10){\vector(-1,0){28}}
\put(60,50){\vector(-1,0){28}}
\put(60,70){\vector(-1,0){28}}
\put(86,10){\vector(-1,0){24}}
\put(86,50){\vector(-1,0){24}}
\put(86,70){\vector(-1,0){24}}
\put(128,10){\vector(-1,0){24}}
\put(128,50){\vector(-1,0){24}}
\put(128,70){\vector(-1,0){24}}
\put(90,10){\circle*{1}}
\put(95,10){\circle*{1}}
\put(100,10){\circle*{1}}
\put(90,50){\circle*{1}}
\put(95,50){\circle*{1}}
\put(100,50){\circle*{1}}
\put(90,70){\circle*{1}}
\put(95,70){\circle*{1}}
\put(100,70){\circle*{1}}
\put(30,25){\circle*{1}}
\put(30,30){\circle*{1}}
\put(30,35){\circle*{1}}
\put(60,25){\circle*{1}}
\put(60,30){\circle*{1}}
\put(60,35){\circle*{1}}
\put(130,25){\circle*{1}}
\put(130,30){\circle*{1}}
\put(130,35){\circle*{1}}
\put(8,30){0}
\put(-23,30){$\infty$}
\put(-7,47){$b_0$}
\put(-7,30){$b_1$}
\put(12,60){$c_{11}$}
\put(12,20){$c_{k1}$}
\put(20,2){$[k,1]$}
\put(24,54){$[2,1]$}
\put(20,74){$[1,1]$}
\put(56,2){$[k,2]$}
\put(56,54){$[2,2]$}
\put(56,74){$[1,2]$}
\put(122,2){$[k,w_k-1]$}
\put(122,54){$[2,w_2-1]$}
\put(122,74){$[1,w_1-1]$}
\put(42,2){$c_{k2}$}
\put(42,54){$c_{22}$}
\put(42,74){$c_{12}$}
\end{picture}
\]
modulo the relations
\begin{equation}\label{e:sqrel}
(\lambda_{i0} b_0 + \lambda_{i1} b_1)c_{i1}=0,
\end{equation}
where $a_i = [\lambda_{i0}:\lambda_{i1}]\in\PP^1$ for $i=1,\dots,k$.
Squids were first studied by Brenner and Butler \cite{BB} in connection with tilting theory.
Recall that a representation $X$ of $S_{D,w}$ is the same as a representation
of the quiver $Q_w$, by vector spaces and linear maps, which satisfies the
relations (\ref{e:sqrel}). The quiver with vertex set $\{\infty,0\}$
and arrows $\{b_0,b_1\}$ is a so-called \emph{Kronecker quiver}, since the
classification of its representations is the same as Kronecker's
classification of pencils of matrices.
We say that a representation of $S_{D,w}$ is \emph{Kronecker-preinjective},
if its restriction to this Kronecker quiver
is preinjective. That is, if $\lambda_0 b_0+\lambda_1 b_1$ is a surjective
linear map for all $[\lambda_0:\lambda_1]\in\PP^1$.

We define the \emph{dimension type} of a representation $X$ of $S_{D,w}$ (or $Q_w$) in a nonstandard way.
If $X_v$ denotes the vector space at vertex $v$, then the dimension type is defined to be the
pair $(\alpha,d)$ where $d = \dim X_\infty$, and $\alpha$ is the dimension vector
for the graph $\Gamma_w$ given by $\alpha_{ij} = \dim X_{[i,j]}$ and $\alpha_0 = \dim X_0 - \dim X_\infty$.

Representations of $Q_w$ of dimension type $(\alpha,d)$ are given by elements
of the set $B \times C$, where
\[
B = M_{d\times(\alpha_0+d)}(K) \times M_{d\times(\alpha_0+d)}(K)
\]
and
\[
C = \prod_{i=1}^k
\left( M_{(\alpha_0+d)\times\alpha_{i1}}(K) \times \prod_{j=2}^{w_i-1} M_{\alpha_{i,j-1}\times\alpha_{ij}}(K) \right).
\]
We denote elements of $B$ and $C$ by $b=(b_0,b_1)$ and $c = (c_{ij})$.
Representations of the squid $S_{D,w}$ are given by elements of the closed subset
\[
S(\alpha,d) = \{ (b,c)\in B\times C \mid \text{$(\lambda_{i0} b_0 + \lambda_{i1} b_1)c_{i1}=0$ for all $i$}\},
\]
and the Kronecker-preinjective representations form a subset $\mathrm{KI}(\alpha,d)$ of $S(\alpha,d)$.
The algebraic group
\[
\GL(\alpha,d) = \GL_d(K) \times \GL_{\alpha_0+d}(K) \times \prod_{i,j} \GL_{\alpha_{ij}}(K)
\]
acts on each of these, with the orbits corresponding to isomorphism classes.

\begin{lem}\label{l:KPirred}
$\mathrm{KI}(\alpha,d)$ is an open subset of $S(\alpha,d)$. It is
an irreducible and smooth variety of dimension
\[
\dim \mathrm{KI}(\alpha,d) = \dim \GL(\alpha,d) - q(\alpha).
\]
\end{lem}

\begin{proof}
Representations of the Kronecker quiver of dimension vector $(d,\alpha_0+d)$
are parametrized by the set $B$.
We show that the preinjective representations form an open subset $U\subseteq B$.
Let $\PP((K^d)^*)$ denote the projective space of nonzero linear forms $\xi$ on $K^d$,
up to rescaling. Let $Z$ be the closed subset of $B\times\PP^1\times \PP((K^d)^*)$
consisting of the elements $((b_0,b_1),[\lambda_0:\lambda_1],\xi)$
with $\xi(\lambda_0 b_0+\lambda_1 b_1)=0$. Clearly $B\setminus U$ is the
image of $Z$ under the projection to $B$, so by projectivity it is closed.
(We thank a referee for this argument, which is shorter than our original one.)

Now $\mathrm{KI}(\alpha,d) = S(\alpha,d) \cap (U\times C)$,
so $\mathrm{KI}(\alpha,d)$ is open in $S(\alpha,d)$.

Let $V = \prod_{i=1}^k M_{d\times \alpha_{i1}}(K)$ and let $h:U\times C\to U\times V$
be the map sending a pair $(b,c)$ to to the element $(b,(\lambda_{i0} b_0 + \lambda_{i1} b_1)c_{i1})$.
Clearly $h$ is a homomorphism of trivial vector bundles over $U$, and
using the fact that $\lambda_{i0} b_0 + \lambda_{i1} b_1$ is a surjective
linear map for all $b\in U$, it is easy to see that
$h$ is surjective on fibres.
It follows that the kernel of $h$ is a vector bundle over $U$ of rank
\[
r = \sum_{i=1}^k \left( \alpha_0\alpha_{i1} + \sum_{j=2}^{w_i-1} \alpha_{i,j-1}\alpha_{ij} \right).
\]
However, this kernel is clearly equal to $\mathrm{KI}(\alpha,d)$. The assertions follow.
\end{proof}

Because $\mathrm{KI}(\alpha,d)$ is irreducible, one can ask about the properties of a general
element of it. Is it indecomposable? What is the dimension of its endomorphism ring?
And so on. Note that the closure of $\mathrm{KI}(\alpha,d)$ is an irreducible component
of $S(\alpha,d)$, so such questions fit into the theory initiated in \cite{CBS}.
Given dimension types $(\beta,e)$ and $(\gamma,f)$, we define $\ext((\beta,e),(\gamma,f))$
to be the general dimension of $\Ext^1(X,Y)$ with $X$ in $\mathrm{KI}(\beta,e)$ and
$Y$ in $\mathrm{KI}(\gamma,f)$.
By the main theorems of \cite{CBS}, one has the following canonical decomposition.

\begin{thm}
Given a dimension type $(\alpha,d)$, there are decompositions
\[
\alpha = \beta^{(1)} + \beta^{(2)} + \dots, \qquad d = d^{(1)}+ d^{(2)}+\dots,
\]
unique up to simultaneous reordering, such that the general element of $\mathrm{KI}(\alpha,d)$
is a direct sum of indecomposables of dimension types $(\beta^{(i)},d^{(i)})$.
These decompositions of $\alpha$ and $d$
are characterized by the property that the general
elements of $\mathrm{KI}(\beta^{(i)},d^{(i)})$ are indecomposable, and
$\ext((\beta^{(i)},d^{(i)}),(\beta^{(j)},d^{(j)}))=0$ for $i\neq j$.
\end{thm}

Following the argument of \cite[Theorem 3.3]{KR}, we have the following result.

\begin{lem}
If $\alpha$ is in the fundamental region, nonisotropic \textup{(}that is $q(\alpha)\neq 0$\textup{)},
or indivisible \textup{(}that is, its components have no common divisor\textup{)},
then the general element of $\mathrm{KI}(\alpha,d)$ is indecomposable.
\end{lem}

\begin{proof}
If the general element of $\mathrm{KI}(\alpha,d)$ decomposes, there are
decompositions $\alpha=\beta+\gamma$ and $d=e+f$, such that the map
\[
\phi:\GL(\alpha,d)\times \mathrm{KI}(\beta,e)\times \mathrm{KI}(\gamma,f)\to \mathrm{KI}(\alpha,d),
\quad
(g,x,y)\mapsto g\cdot(x\oplus y)
\]
is dominant. But $\phi$ is constant on the orbits of the free action
of $\GL(\beta,e)\times \GL(\gamma,f)$ on the domain of $\phi$ given by
\[
(g_1,g_2)\cdot(g,x,y) = (g \begin{pmatrix} g_1 & 0 \\ 0 & g_2 \end{pmatrix}^{-1}, g_1 x, g_2 y).
\]
Thus
\[
\begin{split}
\dim \mathrm{KI}(\alpha,d)\le \left(\dim \GL(\alpha,d)+\dim \mathrm{KI}(\beta,e)+\dim \mathrm{KI}(\gamma,f)\right)
\\
-\left(\dim \GL(\beta,e)+\dim \GL(\gamma,f)\right),
\end{split}
\]
so $q(\alpha)\ge q(\beta)+q(\gamma)$. However, since $\alpha$ is in the fundamental region,
the lemmas in \cite[\S 3.1]{KR} now imply that $\alpha$ is isotropic and divisible.
\end{proof}

Let $A_0$ be the path algebra of the Kronecker quiver with
vertex set $\{\infty,0\}$ and arrows $\{b_0,b_1\}$.
The regular modules for $A_0$ form a stable separating family $\mathcal{T}$,
see \cite[\S 3.2]{Ringel}.
For $1\le i\le k$, let $E_i$ be the regular module of dimension vector (1,1) in which
$\lambda_{i0} b_0 + \lambda_{i1} b_1 = 0$, and let $K_i$ be the full subquiver
of $Q_w$ on the vertices $[i,1],\dots,[i,w_i-1]$. Clearly
$S_{D,w}$ is identified with the tubular extension $A_0[E_i,K_i]_{i=1}^k$ in the sense
of Ringel \cite[\S 4.7]{Ringel}.

Note that since the root of the branch $K_i$ is the unique sink in $K_i$, the length
function $\ell_{K_i}$ (see \cite[\S 4.4]{Ringel}) is equal to $w_i - j$ at vertex $[i,j]$.
This is the same as the dimension vector of the direct sum of all indecomposable
projective representations of $K_i$. Thus, if $M$ is an indecomposable module whose support
is contained in $K_i$, then $\langle \ell_{K_i},\underline{\dim} M\rangle > 0$.

We define module classes $\mathcal{P}_0$, $\mathcal{T}_0$ and $\mathcal{Q}_0$ in
the category $S_{D,w}$-mod as follows.
Let $\mathcal{P}_0$ be the class of preprojective $A_0$-modules,
let $\mathcal{T}_0$ be the class given by all indecomposable modules $M$ such that the underlying Kronecker-module is nonzero
and regular, and let $\mathcal{Q}_0$ be the class of all Kronecker-preinjective modules.

By Ringel \cite[Theorem 4.7(1)]{Ringel}, every indecomposable $S_{D,w}$-module
belongs to $\mathcal{P}_0$, $\mathcal{T}_0$ or $\mathcal{Q}_0$, and
$\mathcal{T}_0$ is a separating tubular family, separating $\mathcal{P}_0$ from $\mathcal{Q}_0$.
In particular, there are no nonzero maps from modules in $\mathcal{Q}_0$ to modules
in $\mathcal{P}_0$ or $\mathcal{T}_0$. Now since all indecomposable
projective modules belong to $\mathcal{P}_0$ or $\mathcal{T}_0$,
Ringel \cite[2.4(1*)]{Ringel} shows that all Kronecker-preinjective
representations have injective dimension $\le 1$.

\begin{lem}\label{l:frinf}
If $\alpha$ is in the fundamental region and $d\ge 0$, then there are infinitely many
isomorphism classes of indecomposable Kronecker-preinjective representations of $S_{D,w}$ of
dimension type $(\alpha,d)$.
\end{lem}

\begin{proof}
First suppose that $\alpha$ is nonisotropic. Thus the general element of
$\mathrm{KI}(\alpha,d)$ is indecomposable. Now isomorphism classes correspond to orbits of
the group $\GL(\alpha,d)$, or actually its quotient $\GL(\alpha,d)/K^*$.
If there were only finitely many isomorphism classes of indecomposables, one would have
\[
\dim \GL(\alpha,d) - q(\alpha) = \dim \mathrm{KI}(\alpha,d) \le \dim \GL(\alpha,d) - 1
\]
so $q(\alpha)\ge 1$, which is impossible for $\alpha$ in the fundamental region.

Now suppose that $\alpha$ is isotropic. By deleting any vertices $v$ of $\Gamma_w$ with
$\alpha_v=0$, we may suppose that $\Gamma_w$ is an extended Dynkin diagram.
Thus, up to reordering, and omitting weights which are equal to 1, we have
that $w$ is of the form (2,2,2,2), (3,3,3), (4,4,2) or (6,3,2).
Thus the squid $A = S_{D,w}$ is a tubular algebra \cite[\S 5]{Ringel}.

Let $\chi_A$ be the quadratic form on $K_0(A)$ for the algebra $A$,
see \cite[\S 2.4]{Ringel}. If $S_v$ denotes the simple $A$-module
corresponding to vertex $v$ in $Q_w$, then
\[
\dim \Ext^1(S_0,S_\infty)=2,
\quad
\dim \Ext^1(S_{[i,1]},S_0)=1,
\]
\[
\dim \Ext^1(S_{[i,j]},S_{[i,j-1]})=1,
\quad
\dim \Ext^2(S_{[i,1]},S_\infty)=1
\]
for $1\le i\le k$ and $1<j\le w_i-1$, and all other Ext spaces between
simples are zero. If $(\alpha,d)$ is a dimension type,
the corresponding element of $K_0(A)$ is
\[
c = d[S_\infty] + (\alpha_0+d)[S_0] + \sum_{i=1}^k\sum_{j=1}^{w_i-1} \alpha_{ij} [S_{[i,j]}],
\]
and using the homological formula for $\chi_A$, one computes
\[
\begin{split}
\chi_A(c)
& = d^2 + (\alpha_0+d)^2 + \sum_{i,j} \alpha_{ij}^2 - 2d(\alpha_0+d) - \sum_i(\alpha_0+d)\alpha_{i1} \\
& \qquad - \sum_{i,j>1} \alpha_{i,j-1}\alpha_{ij} + \sum_i d\alpha_{i1} \\
& = q(\alpha).
\end{split}
\]
Thus if $\alpha$ is isotropic, then $\chi_A(c)=0$, so by \cite[Theorem 5.2.6]{Ringel}
there are infinitely many indecomposable $A$-modules of dimension type $(\alpha,d)$.
\end{proof}

We now consider parabolic bundles $\E$ on $\PP^1$ of weight type $(D,w)$,
where the underlying vector bundle $E$ is algebraic, or equivalently,
if the base field $K$ is the field of complex numbers, holomorphic.

Tilting theory \cite{BB} led first to the theory of derived
equivalences between algebras, and then to derived equivalences between
categories of coherent sheaves on weighted projective lines
and certain algebras, the canonical algebras \cite{Lenzing},
or more generally quasi-tilted algebra \cite{HRS}, including squids.
The following lemma is a special case of the theory.

\begin{lem}\label{l:pbsquidequiv}
There is an equivalence between the category of parabolic bundles
$\E$ on $\PP^1$ of weight type $(D,w)$
such that $E^*$ is generated by global sections, and the category of
Kronecker-preinjective representations of the squid $S_{D,w}$ in which
the linear maps $c_{ij}$ are all injective.

Under this equivalence, parabolic bundles of dimension vector $\alpha$
and degree $-d$ correspond to representations of the squid of
dimension type $(\alpha,d)$.
\end{lem}

\begin{proof}
This follows immediately from the equivalence between the
category of vector bundles $E$ on $\PP^1$ such that $E^*$ is generated by global sections,
and the category of preinjective representations of the Kronecker quiver,
which we now recall.

Let $\OO$ be the trivial bundle on $\PP^1$ and $\OO(1)$
the universal quotient bundle. There is a natural map $\OO^2\to\OO(1)$, and
let $f_0,f_1\in\Hom(\OO,\OO(1))$ be its two components.

The equivalence sends a preinjective representation $X$
of the Kronecker quiver,
\[
X_0
\begin{matrix} \xrightarrow{\;\;\; b_0\;\;\;} \\ \xrightarrow[\;\;\;b_1\;\;\;]{} \end{matrix}
X_\infty,
\]
to the kernel $E$ of the map of vector bundles
\[
f_1\otimes b_0 - f_0\otimes b_1 : \OO\otimes X_0 \to \OO(1)\otimes X_\infty.
\]
The induced map on fibres over $a = [\lambda_0:\lambda_1]$ is the map
\begin{gather*}
X_0 \to [K^2/K(\lambda_0,\lambda_1)] \otimes X_\infty,
\\
x \mapsto [K(\lambda_0,\lambda_1) + (0,1)]\otimes b_0(x) - [K(\lambda_0,\lambda_1) + (1,0)]\otimes b_1(x),
\end{gather*}
which up to a scalar is
the map $\lambda_0 b_0 + \lambda_1 b_1 \in\Hom(X_0,X_\infty)$.
Since $X$ is preinjective, $f_1\otimes b_0 - f_0\otimes b_1$ is surjective on fibres,
and hence $E$ is a vector bundle. Clearly $E^*$ is generated
by global sections. Moreover, using the fact that $\Hom(\OO(1),\OO)=\Ext^1(\OO(1),\OO)=0$,
it is easy to see that this functor is full and faithful. It is also
dense. Namely, recall that for each $i\ge 0$ there is an indecomposable
preinjective Kronecker representation with dimension vector
\[
i+1
\;\;\;
\begin{matrix} \longrightarrow \\ \longrightarrow \end{matrix}
\;\;\;
i
\]
(see for example \cite[\S 3.2]{Ringel}).
Under this functor it gets sent to a vector bundle of rank 1 and
degree $-i$, which must therefore be $\OO(-i)$.

Finally, to see the connection between parabolic bundles and squids,
observe that if $X$ is a Kronecker-preinjective representations of $S_{D,w}$ in which
the linear maps $c_{ij}$ are all injective, and if
$E$ is the corresponding vector bundle, then the squid relation
$(\lambda_{i0} b_0+\lambda_{i1} b_1)c_{i1} = 0$
ensures that the spaces
\[
\Ima(c_{i1}) \supseteq \Ima(c_{i1}c_{i2}) \supseteq \dots
\]
are all contained in the fibre $E_{a_i}$, so define a parabolic bundle $\E$.
Conversely, given a parabolic bundle $\E$ in which the underlying vector
bundle $E$ arises from a preinjective Kronecker module, one can extend this
to a module for $S_{D,w}$ by letting the vector space at $[i,j]$ be $E_{ij}$,
and letting the maps $c_{ij}$ be the natural inclusions.
\end{proof}

A parabolic bundle $\E = (E,E_{ij})$ can be tensored with any vector bundle $F$, giving the
parabolic bundle $\E\otimes F$, with underlying vector bundle $E\otimes F$ and
subspaces $E_{ij}\otimes F_{a_i}$.
Clearly the operation of tensoring with a line bundle is invertible,
so must preserve indecomposability.

\begin{lem}\label{l:infindpbfr}
If $d\in\Z$ and $\alpha$ is in the fundamental region,
then there are infinitely many isomorphism classes of indecomposable
parabolic bundles on $\PP^1$ of weight type $(D,w)$ with dimension vector $\alpha$ and degree $d$.
\end{lem}

\begin{proof}
If $\alpha$ is in the fundamental region, then certainly $\alpha_0>0$, and hence
one can find an integer $N$ such that $d-N\alpha_0 < 0$. By
Lemmas~\ref{l:pbsquidequiv} and \ref{l:frinf} one can find infinitely many
indecomposable parabolic bundles of dimension vector $\alpha$ and degree $d-N\alpha_0$.
Now by tensoring with $\OO(N)$ one obtains
indecomposable parabolic bundles of dimension vector $\alpha$ and degree $d$.
\end{proof}

Note that, instead of using Ringel's work on tubular algebras, it would
have been possible to use the work of Lenzing and Meltzer \cite{LM}.

\section{Riemann-Hilbert correspondence}\label{s:RH}
Let $X$ be a connected Riemann surface and $D=(a_1,\dots,a_k)$ a collection
of distinct points of $X$. We denote by $X\setminus D$ the Riemann surface
obtained by deleting the points of $D$ from $X$.
We write $\OO_X$ for the sheaf of holomorphic functions on $X$, and $\Omega^1_X(\log D)$
for the sheaf of \emph{logarithmic 1-forms} on $X$, that is,
the subsheaf of $j_*\Omega_{X\setminus D}^1$ generated by $\Omega_X^1$ and $dx_i/x_i$,
where $j$ is the inclusion of $X\setminus D$ in $X$ and
$x_i$ is a local coordinate centred at $a_i$.
See \cite[\S II.3]{Deligne}.
If $E$ is a vector bundle on $X$,
a \emph{logarithmic connection}
\[
\nabla:E\to E\otimes \Omega^1_X(\log D)
\]
is a morphism of sheaves of vector spaces which satisfies the Leibnitz rule
\[
\nabla(fe) = e\otimes df + f\nabla(e)
\]
for local sections $f$ of $\OO_X$ and $e$ of $E$. The connection $\nabla$ has residues
\[
\Res_{a_i} \nabla \in \End(E_{a_i}).
\]

Recall that a \emph{transversal} to $\Z$ in $\C$ is a subset $T\subset \C$ with the property
that for any $z\in\C$ there is a unique element $t\in T$ such that $z-t\in\Z$.
Equivalently, $T$ is a transversal if the
function $t\mapsto\exp(-2\pi\sqrt{-1} \, t)$ is a bijection from $T$ to $\C^*$.
An example is given by $T = \{ z\in \C \mid 0 \le \Re e(z) < 1 \}$, but clearly any subset
of $\C$ whose elements never differ by a nonzero integer can be extended to a transversal.

Let $\mathbf{T} = (T_1,\dots,T_k)$ be a collection of transversals.
We say that a logarithmic connection $\nabla:E\to E\otimes \Omega_X^1(\log D)$
\emph{has eigenvalues in $\mathbf{T}$}
if the eigenvalues of $\Res_{a_i}\nabla$ are in $T_i$ for all $i$.
We write $\conn_{D,\mathbf{T}}(X)$ for the category whose objects are the pairs
$(E,\nabla)$ consisting of a vector bundle $E$
on $X$ and a logarithmic connection $\nabla:E\to E\otimes \Omega_X^1(\log D)$ having
residues in $\mathbf{T}$,
and whose morphisms are the vector bundle homomorphisms commuting with the connections.

We need the following version of the Riemann-Hilbert correspondence,
which is presumably already implicit in the work of Deligne \cite{Deligne}
and others.

\begin{thm}\label{t:rh}
Monodromy gives an equivalence from the category $\conn_{D,\mathbf{T}}(X)$
to the category of representations of $\pi_1(X\setminus D)$. Moreover, any
morphism in $\conn_{D,\mathbf{T}}(X)$ has constant rank.
\end{thm}

\begin{proof}
If $M$ is a connected complex manifold, then monodromy gives an
equivalence from the category of vector bundles on $M$
equipped with an integrable holomorphic connection to the category of
representations of $\pi_1(M)$, see \cite[\S I.1,I.2]{Deligne}.
In particular this applies to $X\setminus D$, and since it has complex
dimension 1, integrability is automatic. It thus suffices to prove that
restriction defines an equivalence from $\conn_{D,\mathbf{T}}(X)$ to the
category of vector bundles on $X\setminus D$ equipped with a
holomorphic connection.

Clearly the restriction functor is faithful. To see that it is dense one uses Manin's
construction---a good reference is Malgrange \cite[Theorem 4.4]{Malgrange}.
It thus remains to check that any morphism
\[
\theta: (G'|_{X\setminus D},\nabla'|_{X\setminus D})\to (G|_{X\setminus D},\nabla|_{X\setminus D})
\]
extends to a morphism $(G',\nabla')\to (G,\nabla)$ which has constant rank.
Now this problem is local on $X$, so we may assume that $X$ is a disk, $D$ is its centre
and $G$ and $G'$ are trivial bundles.
Moreover, we may assume that $\nabla$ and $\nabla'$ have connection forms $\omega = \Gamma\frac{dx}{x}$,
$\omega' = \Gamma'\frac{dx}{x}$ given by constant matrices $\Gamma,\Gamma'$
(see \cite{Malgrange}). Then $\theta$ is given by a holomorphic invertible matrix $S$ on $X\setminus D$
which satisfies $dS = S\omega'-\omega S$. Developing $S$ in its Laurent series
$S = \sum_{i=-\infty}^\infty S_i x^i$,
we obtain
\[
\sum_i i S_i x^{i-1} dx = \sum_i (S_i \Gamma' - \Gamma S_i) x^{i-1} dx,
\]
and hence
\[
(\Gamma + i1) S_i = S_i \Gamma'.
\]
By assumption the eigenvalues of $\Gamma$ and $\Gamma'$ belong to the same transversal
to $\Z$ in $\C$, which implies that $S_i=0$ for $i\neq 0$. The result follows.
\end{proof}

In particular Theorem~\ref{t:rh} shows that the category $\conn_{D,\mathbf{T}}(X)$
is abelian. The following fact is well-known.
See for example \cite[Proposition 1.4.1]{Haefliger}.

\begin{lem}\label{l:rhconj}
If $E$ is a vector bundle on $X$ and $\nabla$ is a connection
with eigenvalues in $\mathbf{T}$,
then the monodromy around $a_i$ is conjugate to $\exp(-2\pi \sqrt{-1} \, \Res_{a_i}\nabla)$.
\end{lem}

Now if $X=\PP^1$, we have
\[
\pi_1(\PP^1\setminus D) = G = \langle g_1,\dots,g_k : g_1\dots g_k = 1\rangle
\]
where the $g_i$ are suitable loops from a fixed base point about the $a_i$.
Thus, given any representation $\rho:\pi_1(\PP^1\setminus D)\to\GL_n(\C)$, one
obtains matrices $\rho(g_1),\dots,\rho(g_k)$ whose product is the identity.

\section{Connections on parabolic bundles}
\label{s:Weil}
Let $X$ be a connected Riemann surface, $D=(a_1,\dots,a_k)$ a collection
of distinct points of $X$, $w=(w_1,\dots,w_k)$ a collection of positive integers
and $\zeta\in\C^w$.

If $\E = (E,E_{ij})$ is a parabolic bundle on $X$ of weight type $(D,w)$,
we say that a logarithmic connection
$\nabla : E \to E \otimes \Omega^1_X(\log D)$ is a \emph{$\zeta$-connection}
on $\E$ if
\begin{equation}
\label{parconcond}
(\Res_{a_i} \nabla - \zeta_{ij} 1)(E_{i,j-1}) \subseteq E_{ij}
\end{equation}
for all $1\le i\le k$ and $1\le j\le w_i$,
where by convention $E_{i0}$ is the fibre $E_{a_i}$ and $E_{i,w_i}=0$.
Equivalently, if the $E_{ij}$ are invariant subspaces for $\Res_{a_i} \nabla$,
and $\Res_{a_i} \nabla$ acts on $E_{i,j-1}/E_{ij}$ as multiplication by $\zeta_{ij}$.

If $\nabla : E \to E \otimes \Omega^1_X(\log D)$ is a logarithmic connection
on a vector bundle $E$, then the residues
$\Res_{a_1}\nabla,\dots,\Res_{a_k}\nabla$
form a $k$-tuple of endomorphisms of vector spaces of the same dimension,
and so one can talk about those logarithmic connections whose residues have
type $\zeta$, in the sense of Section~\ref{s:type}, and the associated dimension vector
$\alpha$ with $\alpha_0 = \rank E$.
We write $\conn_{D,w,\zeta}(X)$ for the category whose objects are the pairs
$(E,\nabla)$ consisting of a vector bundle $E$ on $X$ and a logarithmic connection
$\nabla:E\to E\otimes \Omega_X^1(\log D)$ whose residues have type $\zeta$.
Provided that the elements $\zeta_{ij}$ for fixed $i$ never differ
by a nonzero integer, one can find transversals $\mathbf{T} = (T_1,\dots,T_k)$ such
that $\zeta_{ij}\in T_i$ for all $i,j$, and then $\conn_{D,w,\zeta}(X)$ is a full subcategory
of $\conn_{D,\mathbf{T}}(X)$, and is an abelian category.

If $(E,\nabla)$ is an object in $\conn_{D,w,\zeta}(X)$, then $\nabla$ is
a $\zeta$-connection on the parabolic bundle $\E = (E,E_{ij})$, where
\begin{equation}
\label{e:zetaconform}
E_{ij} = \Ima(\Res_{a_i}\nabla - \zeta_{i1}1)\dots(\Res_{a_i}\nabla - \zeta_{ij}1).
\end{equation}
Moreover the dimension vector of $\E$ is the same as the dimension vector of the
residues of $\nabla$.

Conversely, it is clear that if $\nabla$ is a $\zeta$-connection on a
parabolic bundle $\E$, then $(E,\nabla)$ is an object in $\conn_{D,w,\zeta}(X)$.
However, the dimension vector of $\E$ will not be the same as the dimension vector
associated to the residues of $\nabla$ if any of the inclusions in (\ref{parconcond})
are strict.

By analogy with the `parabolic degree' in \cite{MS}, we define the \emph{$\zeta$-degree}
of a parabolic bundle $\E$ of dimension vector $\alpha$ to be
\[
\deg_\zeta \E = \deg E + \zeta*[\alpha].
\]
Recall that the notation $\zeta*[\alpha]$ was defined in the introduction.
The main result of this section is as follows.

\begin{thm}\label{t:parcon}
A parabolic bundle $\E$ on $X$ of weight type $(D,w)$,
has a $\zeta$-connection if and only if
the dimension vector $\alpha'$ of any indecomposable direct summand
$\E'$ of $\E$ satisfies $\deg_\zeta \E'=0$.
\end{thm}

The case when $D$ is empty is a theorem of Weil \cite{Weil} (see also \cite[Theorem 10]{Atiyah}).
The case when the $\zeta_{ij}$ are strictly increasing rational numbers in the interval $[0,1)$
was proved independently by Biswas \cite{Biswas}.
Compare the theorem also with \cite[Theorem 3.3]{CBmm}.

Our proof of Theorem~\ref{t:parcon} relies on a result of Mihai \cite{Mihai1,Mihai2}.
Note that Mihai uses the notion of an `$s$-connection' on a vector bundle $E$
over a complex manifold, where $s$ is a section with simple zeros of a line
bundle $S$. However, taking $s$ to be the natural section of $S = \OO(D)$,
where $D$ now denotes the divisor $\sum_{i=1}^k a_i$, this coincides for a Riemann surface
with the notion of a logarithmic connection.
Bearing in mind that we use the opposite sign convention for the residues,
Theorem 1 of \cite{Mihai1,Mihai2} can
be written for a compact Riemann surface $X$ as follows.

\begin{thm}[Mihai]\label{t:mihai}
Suppose that $E$ is a vector bundle on $X$, and let $\rho_i\in \End(E_{a_i})$
for $1\le i\le k$.
Then there is a logarithmic connection
$\nabla:E\to E\otimes \Omega^1_X(\log D)$ with $\Res_{a_i}(\nabla) = \rho_i$ for all $i$,
if and only if
\[
\sum_{i=1}^k \tr(\rho_i f_{a_i}) = \frac{1}{2\pi\sqrt{-1}} \, \langle b(E),f \rangle
\]
for all $f\in\End(E)$.
\end{thm}

Here $\langle -,-\rangle : \Ext^1(E,E\otimes\Omega^1_X)\times \End(E) \to \C$
is the pairing corresponding to Serre duality, and $b(E)$ is the class of the
Atiyah sequence
\[
0\to E\otimes\Omega^1_X \to B(E)\to E\to 0.
\]
In particular, Atiyah \cite[Proposition 18]{Atiyah} showed that
\[
\langle b(E),1_E \rangle = -2\pi\sqrt{-1} \, \deg E,
\]
and $\langle b(E),f \rangle = 0$ if $f$ is nilpotent.
Moreover, the Corollary to \cite[Proposition 7]{Atiyah} shows
that the construction of $B(E)$ is additive on direct sums,
which implies the following:

\begin{lem}\label{l:atiyahmore}
If $f\in\End(E)$ is the projection onto a direct summand $E'$ of $E$,
then $\langle b(E),f \rangle = -2\pi\sqrt{-1} \, \deg E'$.
\end{lem}

If $\E$ is a parabolic bundle, we write $\iota^\E_{ij}$ for the inclusions $E_{ij}\to E_{i,j-1}$.

\begin{lem}
Given a parabolic bundle $\E$, there is an exact sequence
\[
\bigoplus_{i=1}^k \bigoplus_{j=1}^{w_i-1} \Hom(E_{i,j-1},E_{ij})
\xrightarrow{F}
\End(E)^* \oplus \bigoplus_{i=1}^k \bigoplus_{j=1}^{w_i-1} \End(E_{ij})
\xrightarrow{G}
\End(\E)^* \to 0,
\]
where $F$ sends $\phi=(\phi_{ij})$ with $\phi_{ij}\in\Hom(E_{i,j-1},E_{ij})$
to the element $(\xi,\eta_{ij})$ with $\xi\in\End(E)^*$
and $\eta_{ij}\in \End(E_{ij})$ given by
$\xi(f)=\sum_{i=1}^k \tr(\iota^\E_{i1}\phi_{i1} f_{a_i})$
for $f\in\End(E)$ and
\[
\eta_{ij} = \begin{cases}
\iota^\E_{i,j+1} \phi_{i,j+1} - \phi_{ij} \iota^\E_{ij}
& \text{\textup(if $1\le j< w_i-1$\textup)}
\\
- \phi_{ij} \iota^\E_{ij}
& \text{\textup(if $j=w_i-1$\textup),}
\end{cases}
\]
and $G$ sends $(\xi,\eta_{ij})$ to the map sending $\theta\in\End(\E)$
to
\[
\xi(\theta) + \sum_{i=1}^k \sum_{j=1}^{w_i-1} \tr(\eta_{ij}\theta_{a_i}|_{E_{ij}}).
\]
\end{lem}

\begin{proof}
Clearly there is an exact sequence
\[
0 \to \End(\E)
\xrightarrow{G'}
\End(E) \oplus \bigoplus_{i=1}^k \bigoplus_{j=1}^{w_i-1} \End(E_{ij})
\xrightarrow{F'}
\bigoplus_{i=1}^k \bigoplus_{j=1}^{w_i-1} \Hom(E_{ij},E_{i,j-1})
\]
where $G'$ sends an endomorphism $\theta\in \End(\E)$ to the collection $(\theta, \theta_{a_i}|_{E_{ij}})$,
and $F'$ sends a collection of endomorphisms $(\theta,\theta_{ij})$ to the element whose
component in $\Hom(E_{ij},E_{i,j-1})$ is
$\theta_{i,j-1} \iota^\E_{ij} - \iota^\E_{ij} \theta_{ij}$ for $j\ge 2$,
and
$\theta_{a_i} \iota^\E_{i1} - \iota^\E_{i1} \theta_{i1}$ for $j=1$.
The required sequence is obtained by dualizing this sequence, and using the trace pairing
to identify $\Hom(U,V)^*$ with $\Hom(V,U)$ for any vector spaces $U,V$.
\end{proof}

If $E$ is a vector bundle on $X$, any endomorphism $f\in\End(E)$, induces
an endomorphism $f_a$ of each fibre $E_a$. Clearly the assignment $X\to \C$,
$a\mapsto \tr(f_a)$ is a globally defined holomorphic function on $X$,
hence constant, since $X$ is compact. We denote its common value by $\tr(f)$.

\begin{proof}[of Theorem~\ref{t:parcon}]
Given $\rho_i\in\End(E_{a_i})$, one has
\begin{equation}
\label{e:thicont}
(\rho_i - \zeta_{ij} 1)(E_{i,j-1}) \subseteq E_{ij}
\end{equation}
for all $j$ if and only if there are maps $\phi_{ij}:E_{i,j-1}\to E_{ij}$
for $1\le j\le w_i-1$ satisfying $\rho_i = \iota^\E_{i1} \phi_{i1} + \zeta_{i1}1$
and
\begin{equation}
\label{e:conzetacond}
(\zeta_{ij}-\zeta_{i,j+1}) 1 = \begin{cases}
\iota^\E_{i,j+1} \phi_{i,j+1} - \phi_{ij} \iota^\E_{ij}
& \text{(if $1\le j< w_i-1$)}
\\
- \phi_{ij} \iota^\E_{ij}
& \text{(if $j=w_i-1$).}
\end{cases}
\end{equation}
Namely, if there are $\phi_{ij}$, then an induction on $j$ shows that
\[
(\rho_i - \zeta_{ij}1)|_{E_{i,j-1}} = \phi_{ij},
\]
so it has image contained in $E_{ij}$, and on the other hand,
if $\rho_i$ satisfies (\ref{e:thicont}) one can take $\phi_{ij}$ to be the
map $E_{i,j-1} \to E_{ij}$ induced by $\rho_i - \zeta_{ij} 1$.
Thus by Theorem~\ref{t:mihai}, there is a $\zeta$-connection on $\E$ if and only if there
are maps $\phi_{ij}$ satisfying (\ref{e:conzetacond}) and
\begin{equation}
\label{e:parconrestr}
\sum_{i=1}^k \tr\left( (\iota^E_{i1}\phi_{i1} + \zeta_{i1} 1) f_{a_i}\right) =
\frac{1}{2\pi\sqrt{-1}} \, \langle b(E),f\rangle
\end{equation}
for all $f\in\End(E)$. Now if one defines
\[
(\xi,\eta_{ij})\in \End(E)^* \oplus \bigoplus_{i=1}^k \bigoplus_{j=1}^{w_i-1} \End(E_{ij})
\]
by
\[
\xi(f) = \frac{1}{2\pi\sqrt{-1}} \, \langle b(E),f\rangle - (\sum_{i=1}^k \zeta_{i1}) \, \tr(f),
\]
for $f\in\End(E)$, and $\eta_{ij} = (\zeta_{ij}-\zeta_{i,j+1}) 1$,
then $(\xi,\eta_{ij})\in\Ima(F)$ if and only if there are elements
$\phi_{ij}\in\Hom(E_{i,j-1},E_{ij})$ satisfying
(\ref{e:conzetacond}) and (\ref{e:parconrestr}).

Thus by exactness, there is a $\zeta$-connection on $\E$
if and only if $(\xi,\eta_{ij})\in\Ker(G)$, that is,
\begin{equation}
\label{e:zetaconcond}
\frac{1}{2\pi\sqrt{-1}} \, \langle b(E),\theta\rangle
- (\sum_{i=1}^k \zeta_{i1}) \, \tr(\theta)
+ \sum_{i=1}^k \sum_{j=1}^{w_i-1} (\zeta_{ij}-\zeta_{i,j+1}) \tr( \theta_{a_i}|_{E_{ij}} )
= 0
\end{equation}
for all $\theta\in\End(\E)$.

Suppose that there is a $\zeta$-connection on $\E$,
and let $\theta$ be the projection onto an indecomposable direct
summand $\E' = (E',E'_{ij})$ of $\E$ of dimension vector $\alpha'$. By Lemma~\ref{l:atiyahmore},
equation (\ref{e:zetaconcond}) becomes
\[
-\deg E' - (\sum_{i=1}^k \zeta_{i1}) \, \alpha'_0 +
\sum_{i=1}^k
\sum_{j=1}^{w_i-1} (\zeta_{ij}-\zeta_{i,j+1}) \alpha'_{ij}
= 0.
\]
That is, $\deg E' + \zeta*[\alpha'] = 0$.

For the converse it suffices to show that an indecomposable
parabolic bundle $\E$ of dimension $\alpha$ with $\deg E + \zeta*[\alpha] = 0$
has a $\zeta$-connection. Now $\mathrm{par}_{D,w}(X)$
is an additive category with split idempotents whose Hom spaces are finite-dimensional
vector spaces over $\C$, so any
endomorphism of an indecomposable parabolic bundle can be written as $\theta = \lambda 1+\phi$
with $\lambda\in\C$ and $\phi$ nilpotent. Thus by linearity it suffices to check
equation (\ref{e:zetaconcond}) for $\theta = 1$ and for $\theta$ nilpotent.
The first holds because of the degree condition, the second because all terms are zero.
Thus equation (\ref{e:zetaconcond}) holds for all $\theta$, and hence there is
a $\zeta$-connection on $E$.
\end{proof}

Sometimes we shall use the following special case of Theorem~\ref{t:parcon}.
It also follows directly from \cite[Corollaire 3]{Mihai2}.

\begin{cor}\label{c:parconcor}
If a parabolic bundle $\E$ on $X$
of weight type $(D,w)$ has a $\zeta$-connection,
then $\deg_\zeta \E = 0$.
\end{cor}

\section{Generic eigenvalues}\label{s:generic}
In this section $X=\PP^1$, $D=(a_1,\dots,a_k)$ is a collection
of distinct points of $\PP^1$ and $w=(w_1,\dots,w_k)$ is a collection of positive integers.
Note that in some of the results below $D$ is explicitly given, while in others, such
as Theorem~\ref{t:gencase}, $D$ does not occur in the statement, and it may be chosen
arbitrarily.

Given $\alpha\in\Z^I$,
we say that $\xi\in(\C^*)^w$ is a \emph{generic solution}
to $\xi^{[\alpha]}=1$
provided that for $\beta\in\Z^I$,
\[
\xi^{[\beta]}=1
\Leftrightarrow
\text{$\beta$ is an integer multiple of $\alpha$.}
\]
Given $\alpha\in\Z^I$ and $d\in\Z$,
we say that $\zeta\in\C^w$ is a \emph{generic solution} to $\zeta*[\alpha]=-d$
provided that $\zeta*[\alpha]=-d$ and for $\beta\in\Z^I$,
\[
\zeta*[\beta]\in\Z
\Leftrightarrow
\text{$\beta$ is an integer multiple of $\alpha$.}
\]
In this section we study the existence of representations of $G$ of type $\xi$,
or the existence of connections with residues of type $\zeta$, in these
generic cases. The results will be used in the next section to prove the
main theorems.

\begin{lem}\label{l:gensols}
Suppose that $\alpha\in\Z^I$, $\zeta\in\C^w$ and
$\zeta*[\alpha]=-d\in\Z$. Define $\xi\in(\C^*)^w$ by
$\xi_{ij}=\exp(-2\pi\sqrt{-1}\, \zeta_{ij})$.
Then
$\xi$ is a generic solution to $\xi^{[\alpha]}=1$
if and only if
$\zeta$ is a generic solution to $\zeta*[\alpha]=-d$.
\end{lem}

\begin{proof}
$\xi^{[\beta]} = \exp(-2\pi\sqrt{-1}\, \zeta*[\beta])$, so
$\xi^{[\beta]} = 1\Leftrightarrow\zeta*[\beta]\in\Z$.
\end{proof}

\begin{lem}\label{l:gensolexist}
Suppose $\alpha$ is nonzero.

\textup{\hphantom{i}(i)}~If $\alpha,d$ are coprime,
there is always a generic solution to $\zeta*[\alpha]=-d$.

\textup{(ii)}~There is always a generic solution to $\xi^{[\alpha]}=1$.
\end{lem}

\begin{proof}
For (i), one just needs to remove countably many hyperplanes from
the set $\{ \zeta\in\C^w : \zeta*[\alpha]=-d\}$. For (ii), we use part (i) to
choose a generic solution to $\zeta*[\alpha] = -1$, and then
the previous lemma gives a generic solution to $\xi^{[\alpha]}=1$.
\end{proof}

\begin{thm}\label{t:gencase}
Let $\alpha\in\N^I$ be strict, and suppose
that $\xi$ is a generic solution to $\xi^{[\alpha]}=1$.
Then there is a representation of $G$ of type $\xi$ and dimension $\alpha$
if and only if $\alpha$ is a root. If $\alpha$ is a real root,
this representation is unique up to isomorphism.
\end{thm}

Note that any such representation must be irreducible,
for any subrepresentation also has type $\xi$, and if it has dimension vector
$\beta$ then $\xi^{[\beta]}=1$ by Lemma~\ref{l:repdetcond}.
By genericity, $\beta$ is a multiple of $\alpha$, and in particular the
dimension of the subrepresentation is a multiple of the dimension of the
representation.

\begin{proof}
First note that in case $\alpha$ is a real root, we have already shown the
existence of a representation of $G$ of type $\xi$ and dimension $\alpha$,
unique up to isomorphism, in Theorem~\ref{t:rigidds}.

Suppose that $\alpha$ is a root. We show the existence of a representation of type $\alpha$
by induction.

To get started, this is clear if $\alpha=\epsilon_0$.
Suppose that $\alpha$ is in the fundamental region. Choose a collection of
transversals $\mathbf{T}=(T_1,\dots,T_k)$, and define $\zeta\in\C^w$ with
$\zeta_{ij}\in T_i$ and $\exp(-2\pi\sqrt{-1}\, \zeta_{ij})=\xi_{ij}$.
The condition that $\xi^{[\alpha]}=1$ implies that $d=-\zeta*[\alpha]\in\Z$.
By Lemma~\ref{l:infindpbfr} there is an indecomposable parabolic bundle
$\E$ of dimension vector $\alpha$ and degree $d$.
Now by Theorem~\ref{t:parcon} this parabolic bundle has a $\zeta$-connection
$\nabla$. Clearly the residues of $\nabla$ have type $\zeta$
(but possibly a different dimension vector). By the Riemann-Hilbert
correspondence, Theorem~\ref{t:rh} and Lemma~\ref{l:rhconj}, one obtains a
representation of $G$ of dimension $\alpha_0$ and type $\xi$.
Now by genericity it must have dimension vector $\alpha$.

Now suppose that $\alpha$ is not $\epsilon_0$ or in the fundamental region.
Since it is a root, some reflection $\alpha'=s_v(\alpha)$ is
smaller than $\alpha$. Moreover, we can choose $v$ so that $\alpha'$ is still
strict. Namely, if $\alpha'$ is not strict then $v=0$ and $\alpha'$ is of the form (\ref{e:nonstrictroots}),
so $\alpha$ must be of the form (\ref{e:colldimv}), in which case $s_u(\alpha)<\alpha$
for $u=[\ell,r]$, and $s_u(\alpha)$ is strict.

Observe that $\xi'=r_v'(\xi)$ is
a generic solution to $(\xi')^{[\alpha']}=1$, so by induction there is a
representation $\rho'$ of $G$ of type $\xi'$ and dimension vector $\alpha'$.
By genericity, $\xi^{[\epsilon_v]}\neq 1$. Thus by convolution,
Theorem~\ref{t:conv}, in case $v=0$ or by Lemma~\ref{l:armreflrep} in case $v=[i,j]$,
one obtains a representation $\rho$ of $G$ of type $\xi$ and dimension vector $\alpha$.
(To see that $\rho'$ is $\xi'$-noncollapsing, note that since it is irreducible,
if it weren't, it would be $\xi'$-collapsing. But this is impossible since
$s_v(\alpha')=\alpha>\alpha'$.)

It remains to prove that if there is a representation of $G$ of type $\xi$
and dimension $\alpha$ then $\alpha$ is a positive root. This is a special
case of the next result.
\end{proof}

\begin{thm}
Suppose that $\alpha\in\N^I$ is strict and
that $\xi$ is a generic solution to $\xi^{[\alpha]}=1$.
If there is an indecomposable representation of $G$ of type $\xi$
and dimension $r\alpha$, with $r\ge 1$, then $r\alpha$ is a root.
\end{thm}

\begin{proof}
We prove this by induction on $\alpha$. Let $\rho:G\to \GL(V)$ be the representation.

If $\alpha$ is in the fundamental region, then we're done,
as any multiple of it is a root.
Thus suppose that $\alpha$ is not in the fundamental region.
Since it is strict, it has connected support.
Thus, by the definition of the fundamental region, there is some
vertex $v$ with $\alpha'=s_v(\alpha)<\alpha$.
Clearly $\xi'=r_v'(\xi)$ is a generic solution to $(\xi')^{[\alpha']}=1$.

Suppose one can take $v\neq 0$; say $v=[i,j]$. In this case $\alpha'$ is
strict. Now by Lemma~\ref{l:armreflrep}, $\rho$ is an indecomposable
representation of $G$ of type $\xi'$ and dimension $r\alpha'$. By induction,
$r\alpha'$ is a root, and hence so is $r\alpha$.

Thus we may suppose the $\alpha$ is only reduced by the reflection $s_0$.

If $\alpha$ is a multiple of the coordinate vector $\epsilon_0$,
then $\rho(g_i)=\xi_{i1} 1$ for all $i$, so the indecomposability
of $\rho$ implies that it is 1-dimensional, and hence $r\alpha=\epsilon_0$,
which is a root. Thus we may suppose that $\alpha$ is not a multiple of $\epsilon_0$,
and hence by genericity $\xi^{[\epsilon_0]}\neq 1$.

We show that $\rho$ is $\xi$-noncollapsing. Suppose it has a $\xi$-collapsing
representation $\rho'$ in its top or socle. Now $\rho'$ also has
type $\xi$, and if its dimension vector is $\beta$, then $\xi^{[\beta]}=1$,
so by genericity $\beta$ is a multiple of $\alpha$. Now since $\rho'$ is
$\xi$-collapsing, it is 1-dimensional, so its dimension vector must
actually be equal to $\alpha$. However, the dimension vector of
a $\xi$-collapsing representation is either $\epsilon_0$, which we have
already eliminated, or it is reduced by the reflection at a vertex $v\neq 0$,
which we have also dealt with. Thus $\rho$ is $\xi$-noncollapsing.

Thus by convolution, Theorem~\ref{t:conv}, $\rho$ corresponds to an indecomposable
representation of $G$ of type $\xi'$ and dimension $r\alpha'$.
By induction $r\alpha'$ is a root, and hence so is $r\alpha$.
\end{proof}

These results have the following consequences.

\begin{cor}\label{c:genconnroot}
Suppose that $\alpha\in\N^I$ is strict, $d\in\Z$, and
$\zeta$ is a generic solution to $\zeta*[\alpha]=-d$.
Then there is a vector bundle $E$ of degree $d$ and a logarithmic
connection $\nabla:E\to E\otimes \Omega^1_{\PP^1}(\log D)$
with residues of type $\zeta$ and dimension vector $\alpha$
if and only if $\alpha$ is a root. If $\alpha$ is a
real root, the pair $(E,\nabla)$ is unique up to isomorphism.
\end{cor}

\begin{cor}\label{c:genconnisroot}
Suppose that $\alpha\in\N^I$ is strict, $d\in\Z$, and
$\zeta$ is a generic solution to $\zeta*[\alpha]=-d$.
Let $r\ge 1$.
If there is a vector bundle $E$ of degree $rd$,
a logarithmic connection $\nabla:E\to E\otimes \Omega^1_{\PP^1}(\log D)$
with residues of type $\zeta$ and dimension vector $r\alpha$, and
the pair $(E,\nabla)$ is indecomposable, then $r\alpha$ is a root.
\end{cor}

In both cases one wants to apply the Riemann-Hilbert correspondence. In order to
do so, one needs to find transversals
$\mathbf{T}=(T_1,\dots,T_k)$ such that $\zeta_{ij}\in T_i$ for all $i,j$.
The next lemma shows that this is possible.

\begin{lem}
If $\zeta$ is a generic solution to $\zeta*[\alpha]=-d$, with $\alpha$ strict,
then for fixed $i$, the elements $\zeta_{ij}$ never differ by
an integer.
\end{lem}

\begin{proof}
If $\beta$ is the dimension vector given by (\ref{e:nonstrictroots}),
then $\zeta*[\beta] = \zeta_{\ell,s+1}-\zeta_{\ell r}$.
Now $\beta$ is not strict, so not a multiple of $\alpha$, and hence
this is never an integer by genericity.
\end{proof}

In applying the corollaries, the following lemma will be of use.

\begin{lem}\label{l:genstarconstrict}
Suppose that $\alpha\in\N^I$ is strict, $d\in\Z$, and
$\zeta$ is a generic solution to $\zeta*[\alpha]=-d$.
Let $r\ge 1$.
If $\E$ is a parabolic bundle of dimension vector $r\alpha$
and degree $rd$, and $\nabla$ is a $\zeta$-connection on $\E$,
then
\[
(\Res_{a_i}\nabla - \zeta_{ij}1)(E_{i,j-1}) = E_{ij}
\]
for all $i,j$, so that the residues of $\nabla$ have dimension vector $r\alpha$.
\end{lem}

\begin{proof}
If not, then the residues of $\nabla$ have type $\zeta$, but dimension vector $\beta\neq r\alpha$.
Now equation (\ref{e:zetaconform}) defines a parabolic bundle of dimension vector $\beta$
which has a $\zeta$-connection. Thus by
Corollary~\ref{c:parconcor} we have $\zeta*[\beta]\in\Z$.
Thus by genericity $\beta$ is a multiple of $\alpha$. But since $\beta_0=\rank E = r\alpha_0 \neq 0$,
this implies that $\beta=r\alpha$.
\end{proof}

\section{Proofs of the main theorems}\label{s:mainproofs}
Let $X$, $D$ and $w$ be as in the previous section.

\begin{thm}\label{t:indivparexist}
Suppose that $\alpha$ is a strict root for $\Gamma_w$ and $d\in\Z$.
If $d$ and $\alpha$ are coprime, then there is an
indecomposable parabolic bundle on $\PP^1$ of weight type $(D,w)$ with dimension vector
$\alpha$ and degree $d$.
Moreover, if $\alpha$ is a real root, this indecomposable is
unique up to isomorphism.
\end{thm}

\begin{proof}
By Lemma~\ref{l:gensolexist} we can fix a generic solution $\zeta$ of $\zeta*[\alpha]=-d$.
By Corollary~\ref{c:genconnroot} there is a
vector bundle $E$ of degree $d$ and a logarithmic
connection $\nabla:E\to E\otimes \Omega^1_{\PP^1}(\log D)$
with residues of type $\zeta$ and dimension vector $\alpha$.
By (\ref{e:zetaconform}) this defines a parabolic bundle
$\E$ of dimension vector $\alpha$ equipped with a $\zeta$-connection $\nabla$.
Now $\E$ must be indecomposable by Theorem~\ref{t:parcon} and the genericity of $\zeta$.

Now suppose that $\alpha$ is a real root. By Theorem~\ref{t:parcon}, any indecomposable
parabolic bundle $\mathbf{E}$ of dimension vector $\alpha$ and degree $d$
can be equipped with a $\zeta$-connection.
The residues of $\nabla$ have dimension vector $\alpha$ by Lemma~\ref{l:genstarconstrict}.
The uniqueness thus follows from Corollary~\ref{c:genconnroot}.
\end{proof}

\begin{thm}\label{t:starbunisroot}
If there is an indecomposable parabolic bundle $\E$ on $\PP^1$ of weight type $(D,w)$
with dimension vector $\alpha$ and degree $d$, then $\alpha$ is
a strict root for $\Gamma_w$.
\end{thm}

\begin{proof}
Write $\alpha = r\alpha'$ and $d=rd'$ for some $r\ge 1$ with $\alpha',d'$ coprime.
By Lemma~\ref{l:gensolexist} we can fix a generic solution $\zeta$ of $\zeta*[\alpha']=-d'$.

Since $\mathbf{E}$ is indecomposable and $\zeta*[\alpha]=-d=-\deg E$,
Theorem~\ref{t:parcon} gives the existence of a $\zeta$-connection $\nabla$ on $\E$.
Moreover the residues of $\nabla$ satisfy
\[
(\Res_{a_i}\nabla - \zeta_{ij}1)(E_{i,j-1}) = E_{ij}
\]
for all $i,j$ by Lemma~\ref{l:genstarconstrict},
so the spaces $E_{ij}$ can be reconstructed from $\nabla$, and
hence, since $\E$ is indecomposable, $(E,\nabla)$ is an indecomposable
object of $\mathrm{conn}_{D,w,\zeta}(\PP^1)$.
Thus by Corollary~\ref{c:genconnisroot}, $\alpha$ is a root.
\end{proof}

Together, these last two results prove Theorem~\ref{th:mainparroots}.

\begin{thm}\label{t:conncccexist}
Suppose given conjugacy classes $D_1,\dots,D_k$ of type $\zeta$ and dimension vector $\alpha$.
The following are equivalent

\textup{\hphantom{ii}(i)}~There is a vector bundle $E$ on $\PP^1$ and a logarithmic connection $\nabla:E\to E\otimes \Omega^1_{\PP^1}(\log D)$
such that $\Res_{a_i}\nabla \in \overline{D_i}$

\textup{\hphantom{i}(ii)}~One can write $\alpha$ as a sum of strict roots $\alpha=\beta^{(1)}+\beta^{(2)}+\dots$ with
$\zeta*[\beta^{(p)}]\in\Z$ for all $p$.

\textup{(iii)}~One can write $\alpha$ as a sum of positive roots $\alpha=\beta^{(1)}+\beta^{(2)}+\dots$ with
$\zeta*[\beta^{(p)}]\in\Z$ for all $p$, and $\zeta*[\beta^{(p)}]=0$ for those $p$ with
$\beta^{(p)}$ non-strict.
\end{thm}

\begin{proof}
(i)$\Rightarrow$(ii).
Suppose there are $E$ and $\nabla$ with $\Res_{a_i}\nabla \in \overline{D_i}$.
By the implication (i)$\Rightarrow$(ii) of Theorem~\ref{t:conjclotriple},
$\nabla$ is a $\zeta$-connection on some parabolic bundle $\E$ of dimension vector $\alpha$.
Writing $\E$ as a direct sum of indecomposable parabolic bundles
$\E = \E^{(1)}\oplus \E^{(2)}\oplus \dots$,
one obtains a corresponding decomposition $\alpha=\beta^{(1)}+\beta^{(2)}+\dots$.
Now the $\beta^{(p)}$ are strict roots by the previous theorem, and
$\zeta*[\beta^{(p)}] = -\deg E^{(p)} \in\Z$ by Theorem~\ref{t:parcon}.

(ii)$\Rightarrow$(iii) is vacuous.

(iii)$\Rightarrow$(i).
Suppose that $\alpha$ can be written
as a sum of positive roots $\alpha=\beta^{(1)}+\beta^{(2)}+\dots$ with
$d_p=-\zeta*[\beta^{(p)}]\in\Z$ for all $p$, and $d_p=0$ in case $\beta^{(p)}$ is non-strict.
We may assume that $\beta^{(p)}$ and $d_p$ have no common divisor, for if
$\beta^{(p)}=r\beta'$ and $d_p=rd'$, one can replace $\beta^{(p)}$ by $r$
copies of $\beta'$ in the decomposition of $\alpha$.

Suppose that $\beta^{(1)},\dots,\beta^{(t)}$ are the strict roots involved in
the decomposition, and $\beta^{(t+1)},\dots$ are the non-strict roots.
By Theorem~\ref{t:indivparexist} there are indecomposable parabolic bundles of
dimensions $\beta^{(p)}$ and degrees $d_p$ for $p\le t$,
so by Theorem~\ref{t:parcon} there is a parabolic bundle of dimension
vector $\gamma = \beta^{(1)}+\dots+\beta^{(t)}$ equipped with
a $\zeta$-connection.

For a non-strict root $\beta^{(p)}$ of the form (\ref{e:nonstrictroots}),
the condition $\zeta*[\beta^{(p)}]=0$ says that $\zeta_{\ell r} = \zeta_{\ell,s+1}$.
Thus for any $i$ the sequence
$(\gamma_0,\gamma_{i1},\gamma_{i2},\dots,\gamma_{i,w_i-1})$
is obtained from
$(\alpha_0,\alpha_{i1},\alpha_{i2},\dots,\alpha_{i,w_i-1})$
by a finite number of reductions in the sense of Section~\ref{s:type}, with respect to
various $r,s$ and $(\xi_{i1},\xi_{i2},\dots,\xi_{i,w_i})$.
Thus by Lemma~\ref{l:oddconjclo} the residues of $\nabla$ are in the $\overline{D_i}$.
\end{proof}

Finally we can prove Theorem~\ref{t:clothm}.

\begin{thm}
Suppose given conjugacy classes $C_1,\dots,C_k$ of type $\xi$ and dimension vector $\alpha$.
The following are equivalent.

\textup{\hphantom{ii}(i)}~There is a solution to $A_1\dots A_k = 1$ with $A_i\in\overline{C_i}$.

\textup{\hphantom{i}(ii)}~One can write $\alpha$ as a sum of strict roots $\alpha=\beta^{(1)}+\beta^{(2)}+\dots$ with
$\xi^{[\beta^{(p)}]}=1$ for all $p$.

\textup{(iii)}~One can write $\alpha$ as a sum of positive roots $\alpha=\beta^{(1)}+\beta^{(2)}+\dots$ with
$\xi^{[\beta^{(p)}]}=1$ for all $p$.
\end{thm}

\begin{proof}
Fix a collection of transversals $\mathbf{T}=(T_1,\dots,T_k)$, and define
$\zeta$ by $\xi_{ij}=\exp(-2\pi\sqrt{-1} \, \zeta_{ij})$ and $\zeta_{ij}\in T_i$.
Since $\alpha$ satisfies condition (\ref{e:conjdimveccond}) with respect
to $\xi$, it also satisfies it with respect to $\zeta$, so let
$D_1,\dots,D_k$ be the $k$-tuple of conjugacy classes given by $\zeta,\alpha$.
Clearly $C_i = \exp(-2\pi\sqrt{-1} \, D_i)$.
By monodromy, there is a solution to $A_1\dots A_k = 1$
with $A_i\in\overline{C_i}$ if and only if there is a vector bundle $E$ on $\PP^1$ and a logarithmic
connection $\nabla:E\to E\otimes \Omega^1_{\PP^1}(\log D)$ with $\Res_{a_i}\nabla \in\overline{D_i}$.
Observe that $\xi^{[\beta]}= \exp(-2\pi \sqrt{-1} \, \zeta*[\beta])$,
so $\xi^{[\beta]}=1$ if and only if $\zeta*[\beta]\in\Z$.
Note also that if $\beta$ is a non-strict positive root
and $\zeta*[\beta]\in\Z$ then actually $\zeta*[\beta]=0$.
Namely, if $\beta$ is given by (\ref{e:nonstrictroots}),
then $\zeta*[\beta] = \zeta_{\ell,s+1}-\zeta_{\ell r}$, which cannot
be a nonzero integer since $\zeta_{\ell,s+1}$ and $\zeta_{\ell r}$ belong
to the same transversal.
The result thus follows from Theorem~\ref{t:conncccexist}.
\end{proof}

\end{document}